\documentclass[leqno,11pt]{amsart}

\usepackage{amsmath,amssymb,amsthm}
\usepackage{mathrsfs,amsfonts,functan,extarrows,mathtools}

\newtheorem{lemma}{Lemma}[section]
\newtheorem{theorem}{Theorem}[section]
\newtheorem{definition}{Definition}[section]

\newtheorem{remark}{Remark}[section]

\numberwithin{equation}{section}
\allowdisplaybreaks

\def\la{\langle}
\def\ra{\rangle}
\def\lt{\left}
\def\rt{\right}

\def\no{\nonumber}
\def\ls{\lesssim}

\def\p{\partial}
\def\d{\textrm{\,d}}
\def\deq{\overset{\rm def}{=}}

\def\ep{\varepsilon}
\def\om{v}
\def\Om{\Omega}

\def\gradx{\nabla_{\!x}}
\def\gradv{\nabla_{\!v}}

\Macro{L2}{L^2(M_0 \d \om \d x)}
\Macro{L2om}{L^2(M_0 \d \om)}
\Macro{ify}{L^\infty_x}
\Macro{ify2}{L^\infty_{x,\om}}

\def\nm{\norm*{M_0}}

\newscalprod{Double}{\deltwoarg[#1]{#2}{#3}{\la\!\bigl\la}{\ra\!\bigr\ra_{M_0}}{\,,\,}}
\def\skp{\scalprod{Double}}

\def\nfy{\norm{ify}}
\def\nfyom{\norm{ify2}}

\begin{document}

\title[Self-Organized Model]{Hydrodynamic limits of the kinetic self-organized models}

\author[N. Jiang]{Ning Jiang}
\address[Ning Jiang]{\newline School of Mathematics and Statistics, Wuhan University, Wuhan, 430072, P. R. China}
\email{njiang@whu.edu.cn}

\author[L.-J. Xiong]{Linjie Xiong}
\address[Linjie Xiong]{\newline College of Mathematics and Econometrics, Hunan University, Changsha, 410082, Hunan, P. R. China}
\email{xlj@hnu.edu.cn}

\author[T.-F. Zhang]{Teng-Fei Zhang$^\dag$}
\address[Teng-Fei Zhang]{\newline Yau Mathematical Sciences Center, Tsinghua University, Beijing, 100084, P. R. China}
\email{zhtengf@mail.tsinghua.edu.cn}
\thanks{$^\dag$Corresponding author.}


\begin{abstract}
The self-organized hydrodynamic models can be derived from the kinetic version of the Vicsek model. The formal derivations and local well-posedness of the macroscopic equations are done by Degond and his collaborators. In this paper, we rigorously justify this hydrodynamic limit.

\vspace*{5pt}

\vspace*{5pt}

\end{abstract}

\maketitle
\thispagestyle{empty}



\section{Introduction}

It has been an active research area in recent years on the modeling of interacting agent systems arising in nature, such as bird flocks, fish schools, herds of mammals, etc. They provide fascinating examples of self-organized systems which are able to produce large scale stable coherent structures. Among these models, the Vicsek model \cite{Vicsek-1995} has received particular attention due to its simplicity and the universality of its qualitative features. This model is a discrete particle model which consists of a time-discretized set of ordinary differential equations for the particle position and velocities. The velocities are assumed to  be of constant norm and updated according to an alignment rule, i.e. each agent tries to align its velocity to that of its neighbors in some sensing region. Some angular noise is added to account for stochastic fluctuations. A time-continuous version of this model and its kinetic formulation are available in \cite{Czirok-Vicsek-1996, DM-2008-M3AS}. From the the time-continuous Vicsek model to this kinetic model is rigorously justified in \cite{BCC-2012-aml}.

In \cite{DM-2008-M3AS}, Degond-Motsch derived formally the hydrodynamic limit of the time-continuous Vicsek model through an asymptotic analysis of its kinetic version \cite{BCC-2012-aml}. The resulting model is a system of balanced equations for the density and mean velocity orientation (or polarization vector). This model was later called the Self-Organized Hydrodynamics (SOH). Furthermore, in \cite{DLMP-2013-MAA}, Degond-Liu-Motsch-Panferov derived the modifications of the SOH model by the introduction of the attraction-repulsion force, different scaling assumptions about the size of the sensing region which involve a higher level of non-locality. As proved in \cite{DM-2008-M3AS}, the strict combination of alignment and noise results in the appearance of a pressure term in the SOH model. Introducing an attraction-repulsion force and spanning various scaling assumptions on this force and on the size of the sensing region, they obtain in \cite{DLMP-2013-MAA} a variety of effects that are not encompassed in \cite{DM-2008-M3AS}.

In \cite{DLMP-2013-MAA}, besides the derivation of the macroscopic equations, they prove a local-in-time existence theorem in the 2D case for the viscous model (when the non-local effects are retained) and in the 3D case for the inviscid case (when the non-local effects are omitted). Both proofs are based on a suitable symmetrization of the system and on the energy method.

The main concern of the current paper is on the rigorous justification of the hydrodynamic limit from the Self-Organized Kinetic (SOK) system to the Self-Organized Hydrodynamics (SOH). The main challenge is the lack of conservation laws. To overcome this difficulty, the Generalized Collision Invariants (GCI) are employed to derive the macroscopic equations.   We start from the solution constructed in \cite{DLMP-2013-MAA} and proved that there exists a class of the solutions of SOK  uniformly on a time interval on which the solutions of the SOH are constructed, such that their hydrodynamic limits are the SOH. We employ the methods developed in the fluid limits of the Boltzmann equations, such as \cite{Cafli-1980-CPAM} and \cite{Guo-2006-CPAM}.

The paper is organized as follows. In section 2, we introduce the self-organized models: the time-continuous Vicsek model, the kinetic model and the formal hydrodynamic limits. The main results are stated in section 3. Some mathematical preliminaries are introduced in section 4. In section 5, we estimate the term appear in the expansion, and the main estimates for the remainder equation is presented in section 6. In the last section, building on the estimates in the previous two sections, the proof of the main theorem is completed.

\section{The Self-organized model}
\subsection{Self-propelled particles interacting through alignment}
Let $X_k(t)\in\mathbb{R}^{n}$ and $V_k(t)\in \mathbb{S}^{n-1}$ be the position and velocity of the $k\mbox{-}$th particle at time $t$. The time-continuous version of the Vicsek model is written as follow.
\begin{equation}\label{velocity-E}
  \dot{X}_k(t)= c X_k(t)\,,
\end{equation}
\begin{equation}\label{SDE}
\d V_k(t) = P_{V_k^\perp}\circ(\nu \bar{V}_k(t) \d t + \sqrt{2D} \d B^k_t)\,,
\end{equation}
\begin{equation}\label{average-v}
 \bar{V}_k(t)= \frac{J_k}{|J_k|}\,,\quad J_k=\sum_{j,|X_j-X_k|\leq R}V_j\,,
\end{equation}
where for $V\in \mathbb{S}^{n-1}$, $P_{V^\perp}= Id- V\otimes V$ is the orthogonal projection onto the plane orthogonal to $V$. The equation \eqref{SDE} takes the form of a stochastic differential equation (SDE). The projection operator $P_{V_k^\perp}$ ensures that the resulting solution $V_k(t)$ stays on the unit sphere, provided that the SDE is taken in the Stratonovich sense (which is indicated by the symbol $\circ$). The first term inside the bracket is the interaction. It corresponds to a force acting in the direction $\bar{V}_k$ of intensity $\nu$. The second term is a white noise consisting of independent Brownian motion $B^k_t$ in $\mathbb{R}^2$ of intensity $\sqrt{2D}$. Then, let $N\rightarrow \infty$, a mean field model is obtained. This model is described as follows.
\subsection{Mean-field model and scaling}
The mean-field model describes the evolution of the one-particle distribution function $f(t,x,v)$ at position $x\in\mathbb{R}^n$, with velocity $v\in\mathbb{S}^{n-1}$ at time $t\geq 0$. The model is written as:
\begin{equation}\label{F-P}
\partial_t f + cv\cdot\gradx f = -\gradv\cdot(F_ff)+ D\Delta_{\!v}f\,,
\end{equation}
\begin{equation}\label{force-Ff}
F_f(t,x)=\nu P_{v^\perp}\bar{v}_f(t,x)\,,\quad \bar{v}_f(t,x)=\frac{\mathcal{J}_f(t,x)}{|\mathcal{J}_f(t,x)|}\,,
\end{equation}
\begin{equation}\label{Jf}
\mathcal{J}_f(t,x)= \int_{(y,\omega)\in\mathbb{R}^n\times\mathbb{S}^{n-1}}K(\tfrac{|y-x|}{R})f(t,y,\omega)\omega\,\mathrm{d}y\mathrm{d}\omega\,,
\end{equation}
where the constants $c, \nu, D$ are the same as in \eqref{velocity-E} and \eqref{SDE}.
The equation \eqref{F-P} is a Fokker-Planck equation. The left-hand side expresses the rate of change of $f$ due to the spatial transport of the particle with velocity $cv$ while the first term at the right-hand side denotes the transport in velocity space due to the interaction force $F_f$. The last term at the right-hand side is a velocity diffusion term which arises as a consequence of the Brownian noise in particle velocities. Note that since $v$ lies on the sphere, $\Delta_{\!v}$ stands for the Laplace-Beltrami operator on the sphere. The derivation of the mean-field model \eqref{F-P}, \eqref{force-Ff} and \eqref{Jf} from the discrete system \eqref{velocity-E}, \eqref{SDE} and \eqref{average-v} has been justified in \cite{BCC-2012-aml}.

After nondimensionalization process, (for details, see \cite{DFLMN-2013-Schwartz}), we can write down the scaled self-organized kinetic (SOK) model,
\begin{align}\label{Eq:SOK-0}
  \p_t f + \om \cdot \nabla_x f +\eta_0\ \nabla_\om \cdot \left( P_{\om^\bot} D_{f} f \right)
= \frac{1}{\ep}\left(-\nabla_\om \cdot \left( P_{\om^\bot} \Om_{f} f \right)
 +  d \Delta_\om f\right)\,,
\end{align}
where the parameter $\ep$ denotes the mean free path, i.e., the distance needed by a particle to make a finite change in direction of motion due to the interaction force, and $\eta_0=0,1$ indicate the inviscid and viscous case, respectively, and $d>0$ is a constant. Furthermore, the local current density $j_f$, the local mean direction $\Omega_f$, and the quantity $D_f$ are defined respectively as follows:
\begin{align}
j_f(t,x) \deq & \int_{\mathbb{S}^{n-1}} \om f \d \om, \\
\Omega_f(t,x) \deq & \frac{j_f} {|j_f|} \in \mathbb{S}^{n-1}_x ,\\
D_f(t,x) \deq & k P_{\Omega^\bot_f} \frac{\Delta_x j_f}{|j_f|}, \ k>0.
\end{align}

\subsection{Basic Properties of the SOK Model} 
\label{sub:kinetic_properties_of_the_sok_model}


In this subsection, we list some basic properties of the self-organized kinetic (SOK) model, and we refer the readers to \cite{DFL-2015-ARMA,DFLMN-2013-Schwartz,Frou-2012-M3AS} and references therein for more details.

First we propose a hypothesis: $\Om_f=0$, if $|j_f|=0$. Define the collision operator $\mathcal{Q}$ by
\begin{align}
  \mathcal{Q}(f) = -\nabla_\om \cdot (P_{\om^\bot} \Om_f f) + d \Delta_\om f.
\end{align}
We mention that $\mathcal{Q}(f) =d \Delta_\om f$ in the case $|j_f|=0$ by the above hypothesis.

Next we describe the equilibria of $\mathcal{Q}$, which are expressed by the von Mises-Fisher (VMF) distributions with respect to the local mean orientation $\Om \in \mathbb{S}^{n-1}$, namely,
\begin{align}\no
  \mathcal{E} =& \{ f |\ \mathcal{Q}(f)=0 \}
 = \{\om \mapsto \rho M_\Om (\om),\ \forall \rho \in \mathbb{R}_+,\ \Om \in \mathbb{S}^{n-1} \},
\end{align}
where the VMF distribution is defined as
\begin{align}
  M_\Om (\om) = Z_d^{-1} \exp(\frac{\om \cdot \Om}{d})
\end{align}
with a constant $Z_d=\int_{\om \in \mathbb{S}^{n-1}} \exp(\frac{\om \cdot \Om}{d})\d\om$ independent of $\Om$. The VMF distribution enjoys the following properties:
\begin{itemize}
  \item[i)] $M_\Om (\om)$ is a probability density (i.e., $\int_{\om \in \mathbb{S}^{n-1}} M_\Om (\om)\d\om=1$);
  \item[ii)] The first moment of $M_\Om (\om)$ satisfies
  $$\int_{\om \in \mathbb{S}^{n-1}} \Om M_\Om (\om)\d\om=c_1 \Om,\quad
  c_1(d)= \frac{\int_{\om \in \mathbb{S}^{n-1}} (\om \cdot \Om) \exp(\frac{\om \cdot \Om}{d}) \d\om}{\int_{\om \in \mathbb{S}^{n-1}} \exp(\frac{\om \cdot \Om}{d}) \d\om},
  $$
  where the coefficient $c_1 \in [0,1]$ denotes the order parameter in the study of phase transitions.
\end{itemize}

Note that the formula $P_{\om^\bot} \Om =\nabla_\om (\om \cdot \Om)$ ensures that the collision operator $\mathcal{Q}$ can be rewritten as
\begin{align*}
  \mathcal{Q}(f) = d \nabla_\om \cdot \lt(M_{\Om_f} \nabla_\om \lt(\frac{f}{M_{\Om_f}} \rt)\rt),
\end{align*}
which results in a dissipation relation
\begin{align}
    \int_{\om \in \mathbb{S}^{n-1}} \mathcal{Q}(f) \lt(\frac{f}{M_{\Om_f}} \rt) \d\om
  = -d \int_{\om \in \mathbb{S}^{n-1}} \lt| \nabla_\om \lt(\frac{f}{M_{\Om_f}} \rt) \rt|^2 M_{\Om_f} \d\om
  \le 0.
\end{align}
This implies that $\mathcal{Q}(f)=0$ is equivalent to $f \in \mathcal{E}$.

%
%

One of the main difficulties to derive the macroscopic equations of the SOK model is that it obeys only the conservation law of mass. To recover the missing momentum conservation related to the quantity $\Om(t,x)$, Degond-Motsch introduce the concept of the ``Generalized Collision Invariants'' (GCI) in \cite{DM-2008-M3AS},  see also \cite{DFL-2015-ARMA,Frou-2012-M3AS,DFLMN-2013-Schwartz} its applications.
\begin{definition}[\cite{DM-2008-M3AS}]\label{Lemma-GCI}
  For any given $\Om \in \mathbb{S}^{n-1}$, the linearized collision operator $\mathcal{L}_{\Om}$ is defined as
  \begin{align}
      \mathcal{L}_{\Om} f \deq \Delta_\om f- \frac{1}{d}\nabla_\om \cdot (P_{\om^\bot}\Om f)
    = \nabla_\om \cdot \lt(M_{\Om} \nabla_\om \lt(\frac{f}{M_{\Om}} \rt)\rt).
  \end{align}
  The Generalized Collision Invariants (GCI) are the elements in the null space of $\mathcal{L}_{\Om}$:
  \begin{align}
      \mathcal{N}(\mathcal{L}_{\Om}) \deq & \lt\{ \psi \Big| \int_{\om \in \mathbb{S}^{n-1}} \mathcal{L}_{\Om}f\ \psi \d\om=0,\ \forall f \textrm{ such that } \Om_f= \Om \rt\} \\\no
      =& \lt\{ \psi|\ \mathcal{L}_{\Om}^* \psi(\om) = A\cdot \om \textrm{ with } A\cdot \Om =0 \rt\} \\\no
      =& \lt\{ \om \mapsto h(\om\cdot\Om)A \cdot \om +C \textrm{ with } C \in \mathbb{R},\ A\in \mathbb{R}^n, \textrm{ and } A\cdot \Om =0 \rt\},
  \end{align}
  where the operator $\mathcal{L}_{\Om}^*$ is the adjoint of the linearized operator $\mathcal{L}_{\Om}$, which takes the form,
  \begin{align*}
      \mathcal{L}_{\Om}^* \psi = -\Delta_\om \psi - \frac{1}{d} \Om \cdot \nabla_\om \psi
    = -\frac{1}{M_\Om} \nabla_\om \cdot (M_\Om \nabla_\om \psi).
  \end{align*}
  Here, $h(\om \cdot \Om)=h(\cos \theta)=\frac{g(\theta)}{\sin \theta}$ with $g$ being the unique solution of the elliptic equation $\bar{\mathcal{L}}_\Om^* g(\theta)= \sin\theta$ in the space $V$, where
  \begin{align*}
    & \bar{\mathcal{L}}_\Om^* g(\theta) = - \sin^{2-n}\theta e^{-\frac{\cos\theta}{d}} \frac{\rm d}{\rm d \theta}(\sin^{n-2}\theta e^{\frac{\cos\theta}{d}}g'(\theta)) + \frac{n-2}{\sin^2\theta}g(\theta), \\
    & V = \lt\{ g|\ (n-2) \sin^{\frac{n}{2}-2} \theta g \in L^2(0,\pi),\ \sin^{\frac{n}{2}-1} \theta g \in H^1_0 (0,\pi) \rt\}.
  \end{align*}
\end{definition}
Using the GCI, the macroscopic equations of $\Omega$ can be derived, as stated in the following subsection.


\subsection{Formal derivation of the self-organized hydrodynamics}
\label{sub:hilbert_expansion_for_sok_model}


To study the macroscopic limit of the self-organized kinetic (SOK) model, we rewrite \eqref{Eq:SOK-0} as follows,
\begin{align}\label{eps-SOK}
  \p_t f^\ep + \om \cdot \nabla_x f^\ep +\eta_0 \nabla_\om \cdot \left( P_{\om^\bot} D_{f^\ep} f^\ep \right)
= \frac{1}{\ep} \mathcal{Q}(f^\ep)\,.
\end{align}

We seek for a special class of the solutions of \eqref{eps-SOK} of the form:
\begin{align}\label{f-ansatz}
f^\ep = f_0 + \ep f_1 + \ep^2 f^\ep_2,
\end{align}
with the restriction:
\begin{align}\label{J-restriction}
  j_{f^\ep}= j_{f_0}\,.
\end{align}
As a consequence, $\Om_{f^\ep} =\Om_{f_0} = \Om_0\,.$ Using this $\Omega_0$, we can define the Generalized Collisional Invariants as in Lemma \ref{Lemma-GCI}. Under the restriction \eqref{J-restriction}, the nonlinear equation \eqref{eps-SOK}  becomes {\em linear} as follows.
\begin{align}\label{Eq:Expansion}
  \p_t f^\ep + \om \cdot \nabla_x f^\ep +\eta_0\ \nabla_\om \cdot \left( P_{\om^\bot} D_{0} f^\ep \right)
  = \frac{d}{\ep} \mathcal{L}_{\Om_0} f^\ep\,,
\end{align}
where $D_0 \deq D_{f_0}$. We plug \eqref{f-ansatz} into the equation \eqref{Eq:Expansion}, and collect the same orders, which gives:

\noindent\underline{\bf Order $\mathcal{O}(\frac{1}{\ep}):$} The leading order is
\begin{align}\label{Eq:f-0 order}
  0= d \mathcal{L}_{\Om_0} f_0 = \mathcal{Q}(f_0)\,.
\end{align}
Recalling the equivalence between $\mathcal{Q}(f)=0$ and $f\in \mathcal{E}$, the equation \eqref{Eq:f-0 order} implies
\begin{align}
f_0(t,x,v)= \rho_0(t,x) M_{\Om_0(t,x)}(\om)\,,
\end{align}
for some function $\rho_0=\rho_0(t,x)\,.$ In the rest of the paper, we use the notation $M_0= M_{\Omega_0}\,.$

\noindent\underline{\it Order $\mathcal{O}(1)$} To determine the equations satisfied by the macroscopic variables $(\rho_0, \Om_0)$, we consider the order $O(1)$:
\begin{align}\label{Eq:auxiliary}
  d\mathcal{L}_{\Om_0} f_1
  \deq  \p_t f_0 + \om \cdot \nabla_x f_0 +\eta_0\ \nabla_\om \cdot \left( P_{\om^\bot} D_{0} f_0 \right).
\end{align}

We require the part of $f_1$ in $\mathcal{N}(\mathcal{L}_{\Om_0})$ vanishes, and can solve the part in $\mathcal{N}^\perp(\mathcal{L}_{\Om_0})$ ,
\begin{align}\label{f1-solution}
  f_1 = \frac{1}{d} \mathcal{L}_{\Om_0}^{-1} \left( \p_t f_0 + \om \cdot \nabla_x f_0 +\eta_0\ \nabla_\om \cdot \left( P_{\om^\bot} D_{0} f_0 \right) \right),
\end{align}
under the solvability condition that the right-hand side of \eqref{Eq:auxiliary} lies in $\mathcal{N}^\perp(\mathcal{L}_{\Om_0})$ (for details, see section 5). This gives the following macroscopic equations satisfied by $(\rho_0,\Om_0)$:
\begin{align} \tag{SOH} \label{SOH}
  \left\{
  \begin{array}{l}
\partial_t \rho + c_1 \nabla_x \cdot (\rho \Omega)  = 0, \\[0.5em]
\rho \left( \partial_t \Omega + c_2 \Omega \cdot \nabla_x \Omega \right) + d  \,  P_{\Omega^\bot} \nabla_x \rho = c_3 P_{\Omega^\bot} \Delta_x (\rho \Omega)  ,\\[0.5em]
|\Omega| = 1 ,
  \end{array}
  \right.
\end{align}
with the coefficients
\begin{align*}
  &  c_1 = \frac{\int_0^\pi \cos\theta \exp(\frac{\cos\theta}{d}) \sin^{n-2}\theta \d\theta}{\int_0^\pi \exp(\frac{\cos\theta}{d}) \sin^{n-2}\theta \d\theta}, \quad
    c_2 = \frac{\int_0^\pi \cos\theta h(\cos\theta) \exp(\frac{\cos\theta}{d}) \sin^n\theta \d\theta}{\int_0^\pi h(\cos\theta) \exp(\frac{\cos\theta}{d}) \sin^n\theta \d\theta},\\
  & c_3 = \eta_0 k[(n-1)d+ c_2].
\end{align*}
This is the hydrodynamic model, which we call self-organized hydrodynamic (SOH) system. We refer to \cite{DFL-2015-ARMA,DFLMN-2013-Schwartz,Frou-2012-M3AS} for the derivation of the SOH system and omit the details here. Now the equation of the remainder $f^\ep_2$ is
\begin{equation}\label{remainder-equation}
\begin{aligned}
&\p_t f_2^\ep + \om \cdot \nabla_x f_2^\ep +\eta_0\ \nabla_\om \cdot \left( P_{\om^\bot} D_{0} f_2^\ep \right) - \frac{d}{\ep} \mathcal{L}_{\Om_0} f^\ep_2 \\
=& -\left\{\p_t f_1 + \om \cdot \nabla_x f_1 +\eta_0\ \nabla_\om \cdot \left( P_{\om^\bot} D_{0} f_1 \right)\right\}\,.
\end{aligned}
\end{equation}
For notational simplicity, we set $d \equiv 1$ in the rest of the paper. By setting $f_2^\ep = \widetilde{f}_2^\ep M_0$, the equation \eqref{remainder-equation} is reduced to the equation of $\widetilde{f}_2^\ep$:
\begin{align}\label{Remainder-Equation-1}
  \p_t \widetilde{f}_2^\ep + \om \cdot \nabla_x \widetilde{f}_2^\ep +\eta_0\, \nabla_\om \cdot \left( P_{\om^\bot} D_{0} \widetilde{f}_2^\ep \right)
  + \frac{1}{\ep} \mathcal{L}_0 \widetilde{f}_2^\ep
  = h_0 \widetilde{f}_2^\ep + \frac{1}{\ep} h_1\,,
\end{align}
where
\begin{align}
  \mathcal{L}_0 f =-\frac{1}{M_0}\mathcal{L}_{\Omega_0}(M_0 f)= - \frac{1}{M_0} \nabla_\om \cdot (M_0 \nabla_\om f)\,,
\end{align}
\begin{equation}
  h_1 =  - \frac{1}{M_0} \left[ \p_t f_1 + \om \cdot \nabla_x f_1 +\eta_0\, \nabla_\om \cdot \left( P_{\om^\bot} D_{0} f_1 \right)  \right]\,,
\end{equation}
and
\begin{equation}
\begin{aligned}
  h_0 = & - \frac{1}{M_0} \left[ \p_t M_0 + \om \cdot \nabla_x M_0 +\eta_0\, \nabla_\om \cdot \left( P_{\om^\bot} D_{0} M_0 \right)  \right]\\
  =&- [ \om \cdot \p_t \Om_0 + \om \otimes \om: \nabla_x \Om + \eta_0 D_0 \cdot P_{\om^\bot} \Om_0]\,.
\end{aligned}
\end{equation}
In the rest of the paper, we work on the remainder equation \eqref{Remainder-Equation-1}.

\section{Main results} 
\label{sec:main_results}
In this section, we state our main result. We first introduce the existence result of Degond-Liu-Motsch-Panferov, on which our result is built. First, we introduce the Cauchy problem of self-organized hydrodynamic (SOH) system:
\begin{align} \tag{$\rm SOH$} \label{SOH-Cauchy}
  \left\{
  \begin{array}{l}
\partial_t \rho + c_1 \nabla_x \cdot (\rho \Omega)  = 0,\ x \in \mathbb{T},\ t >0, \\[2pt]
\rho \left( \partial_t \Omega + c_2 \Omega \cdot \nabla_x \Omega \right) + d\,  P_{\Omega^\bot} \nabla_x \rho = c_3 P_{\Omega^\bot} \Delta_x (\rho \Omega) ,  \\[2pt]
|\Omega| = 1, \\[2pt]
\rho|_{t=0}= \rho^{in}\ge c_0>0,\ \Om|_{t=0} =\Om^{in}, |\Omega^{in}| = 1.
  \end{array}
  \right.
\end{align}
Here $\mathbb{T}$ denotes $\mathbb{T}^2$ or $\mathbb{T}^3$.

  On the one hand, the SOH system evidently bears many similarities with the isentropic compressible Navier-Stokes (NS) system. And on the other hand, it also has some different properties. The first important difference is that the SOH system obeys the geometric constraint $|\Om|=1$ which requires the velocity $\Om$ to be of unit norm. The second important difference is that, generally speaking, the coefficients $c_i\ (i=1,2)$ are different. Indeed, we have $0<c_2 \le c_1$ (the equality holds iff $d=0$), see \cite{DFLMN-2013-Schwartz,DM-2008-M3AS}.

The existence and uniqueness of the solution to the Cauchy problem of SOH system have been established in \cite{DLMP-2013-MAA}, we quote the results as follows,
\begin{theorem}[\cite{DLMP-2013-MAA}]\label{DLMP-MAA}
 Let $\mathbb{T}=\mathbb{T}^{n}\,,$ given $\rho^{in}\geq c_0 >0$, $\Omega^{in}\in \mathbb{S}^{n-1}$, i.e.
       \begin{itemize}
          \item for $n=2$, $\Omega^{in}=(\cos\varphi^{in}\,, \sin\varphi^{in})\,,  \varphi^{in}\in [0,2\pi]\,;$
          \item for $n=3$, $\Omega^{in}=(\sin\theta^{in}\cos\varphi^{in}\,, \sin\theta^{in}\sin\varphi^{in}\,,\cos\theta^{in})$, $\theta^{in}\in [0,\pi]\,, \varphi^{in}\in [0,2\pi]$. Furthermore, assume $ \sin\theta^{in}>0\,,$ and $c_3=0\,,$ (i.e. $\eta_0=0\,.$)
       \end{itemize}
 with $(\rho^{in}\,,\varphi^{in}\,,\theta^{in})\in H^m\,, m > \frac{n}{2}+1\,.$ Then, there exists $T>0$, such that the Cauchy problem of SOH system with initial data $(\rho^{in}\,,\Omega^{in})$ has a unique solution $(\rho_0\,,\Omega_0)\in L^\infty([0,T], H^m(\mathbb{T}))\cap H^1([0,T],H^{m-1}(\mathbb{T}))$, and $\rho_0 >0\,.$
\end{theorem}

Our hydrodynamic limit result builds on the above theorem, so make the same assumptions on the initial data, namely,

{\noindent \underline{Assumption (A)}:}
\begin{enumerate}
  \item In 2D case, we consider $c_3 \ge 0$ including both viscous and inviscid cases. The initial data $(\rho^{in}, \Om^{in}) \in H^m(\mathbb{T})$ are smooth enough as required. Besides, $\rho^{in}$ has a positive low bound.
  \item In 3D case, we consider only the inviscid case $c_3=0$, which arises from the coefficient $\eta_0=0$. (Indeed, $c_3= \eta_0 \widetilde{c}_3$). Besides, $\rho^{in}$ has a positive low bound and $\Om^{in} \neq (0,0,1)$ (corresponding to the above constraint $\sin \theta^{in}>0$).
\end{enumerate}

Now we state the main result of this paper on the hydrodynamic limit from the self-organized kinetic (SOK) equation to the self-organized hydrodynamic (SOH).
\begin{theorem}\label{Thm: main}
  Let $(\rho_0, \Om_0) \in L^\infty([0,T], H^m(\mathbb{T}))\cap H^1([0,T],H^{m-1}(\mathbb{T}))\ (m>13)$ be the solutions provided by Theorem \ref{DLMP-MAA} to the Cauchy problem of the self-organized hydrodynamic system \eqref{SOH-Cauchy} with initial datum $(\rho^{in}, \Om^{in})$ satisfying the assumption (A). Let $f_0=\rho_0M_{\Omega_0}$, and $f_1(t,x,v)\in\mathcal{N}^\perp(\mathcal{L}_{\Omega_0})$ be the unique solution given in \eqref{f1-solution}.

  Furthermore, assume $f^{\ep,in}(x,\om)= \rho^{in}(x) M_{\Om^{in}(x)}(\om) + \ep f_1(0,x,\om) + \ep^2 f_2^{\ep,in}(x,\om)$ with the bound $\|f_2^{\ep,in}(x,\om)\|_{H^2_{x,\om}} \le C$.

  Then there exists $\ep_0 >0$ such that, for any $\ep \in (0,\ep_0)$, the self-organized kinetic equation \eqref{eps-SOK} admits a unique solution $f^\ep(t,x,\om)\in C([0,T]; H^2(\mathbb{T} \times \mathbb{S}^{n-1}))$ in the class $j_{f^\ep}= j_{f_0}$, of the form
  \begin{align}
    f^\ep(t,x,\om)= \rho_0(t,x) M_{\Om_0(t,x)}(\om) + \ep f_1(t,x,\om) + \ep^2 f_2^\ep(t,x,\om),
  \end{align}
  with $f_2^\ep(t,x,\om)$ satisfying
  \begin{align}\label{Eq: conclusion}
    \ep^\frac{1}{2} \|f_2^\ep\|_{L^2_{x,\om}} + \ep \|f_2^\ep\|_{H^1_{x,\om}}+ \ep^\frac{3}{2} \|f_2^\ep\|_{H^2_{x,\om}} \le C,
  \end{align}
  where the constant $C$ is independent of $\ep \in (0,\ep_0) $ and $t\in[0,T] $.
\end{theorem}

\begin{remark}
  Indeed, we can get a more generic result about the higher order diffusion expansion,
  \begin{align}
    f^\ep = \rho_0 M_{\Om_0} + \ep f_1 + \ep^2 f_2+ \cdots +\ep^{n-1} f_{n-1} + \ep^n f_n^\ep,
  \end{align}
  if it holds initially with the bound $\|f_n^{\ep,in}\|_{H^n_{x,\om}} \le C$. For $i=1,2,\cdots, n-1$, $f_i$ are determined by the equation
  \begin{align*}
    \mathcal{L}_{\Om_0} f_i = \p_t f_{i-1} + \om \cdot \nabla_x f_{i-1} + \eta_0 \nabla_\om \cdot (P_{\om^\bot} D_0 f_{i-1}),
  \end{align*}
  and the $n$-order remainder $f_n^\ep$ satisfies
  \begin{align}
    \sum_{0\le s \le n} \ep^\frac{s+1}{2} \|f_n^\ep\|_{H^s_{x,\om}} \le C,
  \end{align}
  with the constant $C$ independent of $\ep \in (0,\ep_0) $ and $t\in[0,T]$.

  Proving this high order result works exactly the same as for Theorem \ref{Thm: main}, and relies heavily on a high order a priori estimate \eqref{Esm:a priori-high order}, see the discussions in Remarks \ref{Rema:a priori-high order} and \ref{Rema:auxiliary-high order} below.
\end{remark}


\section{Preliminaries} 
\label{sec:preliminaries}



In the following context, we will do estimates in the weighted spaces $L^2(M_0\d \om \d x) $. For the sake of simplicity, we use the notations $\| \cdot \|$ and $\| \cdot \|_{M_0} $ to denote the norms $\| \cdot \|_{L^2(\d \om \d x)} $ and $\| \cdot \|_{L^2(M_0\d \om \d x)} $, respectively. Note that the two norms are equivalent since $M_0$ is bounded from up and below. Similarly, the notations $\la\!\la \cdot, \cdot \ra\!\ra$ and $\skp{\cdot}{\cdot} $ are adapted to stand for the inner products $\la \cdot, \cdot \ra_{L^2(\d \om \d x)}$ and $\la \cdot, \cdot \ra_{L^2(M_0\d \om \d x)}$, respectively.

We will use frequently the following facts on the sphere, for a constant vector $V\in \mathbb{R}^n$,
\begin{align}
  \nabla_\om (\om \cdot V) = (Id - \om \otimes \om) V = P_{\om^\bot} V, \\
  \nabla_\om \cdot ((Id - \om \otimes \om) V) =-(n-1) \om \cdot V,
\end{align}
where $\nabla_\om$ and $\nabla_\om \cdot$ are used to denote the tangential gradient (and, divergence) operator on the sphere. Moreover, we have some other useful formulas,
\begin{align}
  \int \nabla_\om f \d\om = -(n-1) \int \om f \d\om, \\
  \int \om \nabla_\om \cdot A(\om) \d\om = - \int A(\om) \d\om,
\end{align}
with $A(\om) $ any smooth tangent vector field.


Now define
\begin{align*}
  \la\la g \ra\ra= &\iint g \d\om \d x,\qquad \la g \ra = \int g \d\om; \\\no
  \la\la g \ra\ra_{M_0} = &\iint g M_0 \d\om \d x,\quad \la g \ra_{M_0} = \int g M_0 \d\om.
\end{align*}

We give a lemma about the Poincar\'e inequality on the sphere.
\begin{lemma}[Appendix of (\cite{DFL-2013-JNS})] We have the following weighted Poincar\'e inequality, for $g \in H^1(\mathbb{S}^{n-1})$,
\begin{align}
  \la |\nabla_\om g|^2 \ra_{M_0} \ge \Lambda \la (g- \la g \ra_{M_0})^2 \ra_{M_0},
\end{align}
where $\Lambda$ is the Poincar\'e constant independent of $\Om_0$.
\end{lemma}

Due to the compatibility condition of the Hilbert expansion $j_{f^\ep} = j_{f_0},\ \Om_{f^\ep} =\Om_{f_0} $, we have
\begin{align}
  \la \widetilde{f}_1 \ra_{M_0} = \la f_1 \ra =0,\ \la \om \widetilde{f_1} \ra_{M_0} = \la \om f_1 \ra =0; \\
  \la \widetilde{f}_2^\ep \ra_{M_0} = \la f_2^\ep \ra =0,\ \la \om \widetilde{f}_2^\ep \ra_{M_0} = \la \om f_2^\ep \ra =0.
\end{align}
Hence it is convenient to introduce the mean free spaces $\dot{L}^2_{M_0} (\mathbb{S}^{n-1}) \subset L^2 (\mathbb{S}^{n-1})$ composed of functions $g$ satisfying $\la g \ra_{M_0} =0$.

As pointed out in \cite{DFL-2015-ARMA}, the linearized operator $\mathcal{L}_0$ is a self-adjoint operator under the scalar product $\scalprod{L2om}{g_1}{g_2} = \int g_1 g_2 M_0 \d\om$, since we have
\begin{align*}
\scalprod{L2om}{g_1}{\mathcal{L}_0 g_2} = \scalprod{L2om}{\nabla_\om g_1}{\nabla_\om g_2}.
\end{align*}

We quote the definition of an operator $\mathcal{L}_0^s$ and an equivalent Sobolev norm on the sphere by spectral decomposition, as in \cite{DFL-2015-ARMA},
\begin{align}
  \|g\|^2_{\dot H^s_{M_0}(\mathbb{S}^{n-1})} \deq  \scalprod{L2om}{g}{\mathcal{L}_0^s g} .
\end{align}
Then we have the following lemma.
\begin{lemma} [Lemma 2 of \cite{DFL-2015-ARMA}] \label{lemm: Liu-ARMA}
  For $g \in \dot H^s_{M_0}(\mathbb{S}^{n-1})= H^s(\mathbb{S}^{n-1}) \cap \dot L^2(\mathbb{S}^{n-1})$ and $s\ge 0$, we have
  \begin{align*}
   \|g\|^2_{\dot H^s_{M_0}(\mathbb{S}^{n-1})} \sim \|g\|^2_{H^s(\mathbb{S}^{n-1})}.
  \end{align*}
Here, we use the notation $A \sim B$ to denote there exists a universal constant $C>0$, such that $C^{-1} B \le A \le C B$.

  Furthermore, given $g \in \dot H^s_{M_0}(\mathbb{S}^{n-1})$, there exists a constant $C$ such that,
  \begin{align}\label{eq:high order commutator}
    |\scalprod{L2om}{\mathcal{L}_0^s g}{\nabla_\om g} | \le C \|g\|^2_{\dot H^s_{M_0}(\mathbb{S}^{n-1})}.
  \end{align}
\end{lemma}


\section{Estimates of $f_1$}
\label{sec:Estimate_for_auxiliary}

Recalling the equation \eqref{Eq:auxiliary} of $f_1$,
\begin{align*}
  \mathcal{L}_{\Om_0} f_1
  = \p_t f_0 + \om \cdot \nabla_x f_0 +\eta_0\ \nabla_\om \cdot \left( P_{\om^\bot} D_{0} f_0 \right)\,,
\end{align*}
we introduce a new function $\widetilde{f}_1$ satisfying $f_1 = \widetilde{f}_1 M_0$, then the formula
$$\frac{1}{M_0} \mathcal{L}_{\Om_0} f_1 = \frac{1}{M_0} \nabla_\om \cdot (M_0 \nabla_\om \frac{f_1}{M_0})
= \frac{1}{M_0} \nabla_\om \cdot (M_0 \nabla_\om \widetilde{f}_1) = - \mathcal{L}_0 \widetilde{f}_1
$$
gives the equation (the symbol $\sim$ has been dropped for brevity):
\begin{align}\label{Eq:auxiliary-tilde}
  \mathcal{L}_0 f_1 = - \frac{1}{M_0} [\p_t f_0 + \om \cdot \nabla_x f_0 +\eta_0\ \nabla_\om \cdot \left( P_{\om^\bot} D_{0} f_0 \right)].
\end{align}

From the self-adjoint property of $\mathcal{L}_0$, the existence and uniqueness of $f_1$ can be easily established by combining the Lax-Milgram theorem and the Poincar\'e inequality, we omit the details. Next we give a lemma stating the boundedness of the function $f_1$.


\begin{lemma}\label{lemm: Estimate for auxiliary}
  There exists a constant $C>0$ depending on $\|\rho_0\|_{L^\infty(0,T;H^m(\mathbb{T}))} $ and $\|\Om_0\|_{L^\infty(0,T;H^m(\mathbb{T}))} $ with $m>13$, such that for any $t\in[0,T]$,
  \begin{align}
    \|f_1\|_{H^3_{M_0}}(t) + \|\p_t f_1\|_{H^3_{M_0}}(t) \le C.
  \end{align}
\end{lemma}

\proof We split the proof into four steps.

\noindent\textbf{Step I}: Estimation for $\nm{\nabla_\om f_1}$ (hencely for $\nm{f_1}$ by the Poincar\'e inequality on the sphere).

To deal with the first two terms on the right hand side of the equation \eqref{Eq:auxiliary-tilde}, it follows from the Poincar\'e inequality on the sphere that,
\begin{align*}
  & \skp{M_0^{-1} (\p_t f_0 + \om \cdot \nabla_x f_0)}{f_1}
 = \la \p_t f_0 + \om \cdot \nabla_x f_0, f_1 \ra_{L^2(\d\om \d x)} \\
 \le & \|\p_t f_0 + \om \cdot \nabla_x f_0\|\ \|f_1\|
 \ls \|\p_t f_0 + \om \cdot \nabla_x f_0\|\ \|f_1\|_{M_0} \\
 \le & \frac{1}{4} \nm{\nabla_\om f_1}^2 + C\|\p_t f_0 + \om \cdot \nabla_x f_0\|^2.
\end{align*}

As for the third term, we have
\begin{align*}
   & \eta_0 \skp{M_0^{-1} \nabla_\om \cdot (P_{\om^\bot} D_0 f_0)}{f_1}
 \le \eta_0 \nm{M_0^{-1} \nabla_\om \cdot (P_{\om^\bot} D_0 f_0)} \nm{f_1} \\
 \ls & \eta_0 \nm{\nabla_\om \cdot (P_{\om^\bot} D_0 f_0)} \nm{\nabla_\om f_1}\\
 \le & \frac{1}{4} \nm{\nabla_\om f_1}^2 + \eta_0 C \|\nabla_\om \cdot (P_{\om^\bot} D_0 f_0)\|^2.
\end{align*}

Recall that the linearized operator $\mathcal{L}_0$ is nonnegative, i.e.,
\begin{align*}
  \skp{\mathcal{L}_0 f_1}{f_1} = \nm{\nabla_\om f_1}^2,
\end{align*}
then combining the above three estimates, we get
\begin{align}
  \frac{1}{2} \nm{\nabla_\om f_1}^2
 \ls \|\p_t f_0 + \om \cdot \nabla_x f_0\|^2 + \eta_0 \|\nabla_\om \cdot (P_{\om^\bot} D_0 f_0)\|^2.
\end{align}

The fact $f_0 = \rho_0 M_0$ yields that
\begin{align*}
  \p_t f_0=& M_0 \p_t \rho_0 + \rho_0 M_0 \om \cdot \p_t \Om_0, \\
  \nabla_x f_0 =& M_0 \nabla_x \rho_0 + \rho_0 M_0 \om \cdot \nabla_x \Om_0,
\end{align*}
together with the assumption for the (SOH) model $(\rho_0, \Om_0) \in L^\infty(0,T; H^m(\mathbb{T}))$ with $m>4$, the above two equalities result in
\begin{align*}
  \|\p_t f_0 + \om \cdot \nabla_x f_0\|^2 \le C,
\end{align*}
where the constant $C$ depends on the value of $\|\rho_0\|_{L^\infty(0,T;H^m(\mathbb{T}))} $ and $\|\Om_0\|_{L^\infty(0,T;H^m(\mathbb{T}))} $.

Similar argument applied to the term $\nabla_\om \cdot (P_{\om^\bot} D_0 f_0) = -(n-1) D_0 \om f_0 + D_0 \nabla_\om f_0$ shows
\begin{align*}
  \eta_0 \|\nabla_\om \cdot (P_{\om^\bot} D_0 f_0)\|^2 \le C.
\end{align*}

Thus we are lead to the conclusion
\begin{align}
  \nm{f_1}^2  \ls \nm{\nabla_\om f_1}^2 \le C.
\end{align}

\noindent\textbf{Step II}: Applying the operator $\nabla_x$ to the equation \eqref{Eq:auxiliary-tilde}, and taking $\m{L2}$ scalar product with $\nabla_x f_1$, we can get the estimation for $\nm{\nabla_\om \nabla_x f_1} $ (and hence, for $\nm{\nabla_x f_1} $). We only collect the controls for both sides of the equation \eqref{Eq:auxiliary-tilde} as follows,
\begin{align*}
  & \skp{\nabla_x \mathcal{L}_0 f_1}{\nabla_x f_1} \\
  =& \skp{\mathcal{L}_0 \nabla_x f_1}{\nabla_x f_1} - \skp{\nabla_x \Om_0 \nabla_\om f_1}{\nabla_x f_1} \\
  \ge & \nm{\nabla_\om \nabla_x f_1}^2 - \nfy{\nabla_x \Om_0} \nm{\nabla_\om f_1} \nm{\nabla_x f_1}\\
  \ge & \frac{3}{4} \nm{\nabla_\om \nabla_x f_1}^2 - C \nfy{\nabla_x \Om_0}^2 \nm{\nabla_\om f_1}^2,
\\
  & \skp{\nabla_x \{ M_0^{-1} [\p_t f_0 + \om \cdot \nabla_x f_0 + \eta_0\ \nabla_\om \cdot \left( P_{\om^\bot} D_{0} f_0 \right)] \} }{\nabla_x f_1} \\
  \le & \nm{\nabla_x \{ M_0^{-1} [\p_t f_0 + \om \cdot \nabla_x f_0 + \eta_0\ \nabla_\om \cdot \left( P_{\om^\bot} D_{0} f_0 \right)] \}} \nm{\nabla_x f_1} \\
  \le & \frac{1}{4} \nm{\nabla_\om \nabla_x f_1}^2 + C \nm{\nabla_x \{ M_0^{-1} [\p_t f_0 + \om \cdot \nabla_x f_0 + \eta_0 \nabla_\om \cdot \left( P_{\om^\bot} D_{0} f_0 \right)] \}}^2.
\end{align*}

The above two inequalities together implies
\begin{align}
  \frac{1}{2} \nm{\nabla_\om \nabla_x f_1}^2 \ls & \nfy{\nabla_x \Om_0}^2 \nm{\nabla_\om f_1}^2 \\\no
  &+ \nm{\nabla_x \{ M_0^{-1} [\p_t f_0 + \om \cdot \nabla_x f_0 + \eta_0 \nabla_\om \cdot \left( P_{\om^\bot} D_{0} f_0 \right)] \}}^2,
\end{align}
combining with the fact $\nm{\nabla_\om f_1} \le C$ from step I and the assumption of $(\rho_0, \Om_0) \in L^\infty(0,T; H^m(\mathbb{T}))$ with $m>4$, this gives the result
\begin{align}
  \nm{\nabla_x f_1}^2 \ls \nm{\nabla_\om \nabla_x f_1}^2 \le C.
\end{align}

\noindent\textbf{Step III}: The estimations for $\nm{\nabla_\om \p_t f_1} $ (and $\nm{\p_t f_1} $) is similar to that of step II, and we finally get
\begin{align}
  & \nm{\p_t f_1}^2 \ls \nm{\nabla_\om \p_t f_1}^2 \\\no
 \ls & \nfy{\p_t \Om_0}^2 \nm{\nabla_\om f_1}^2
  + \nm{\p_t \{ M_0^{-1} [\p_t f_0 + \om \cdot \nabla_x f_0 + \eta_0 \nabla_\om \cdot \left( P_{\om^\bot} D_{0} f_0 \right)] \}}^2 \\\no
 \ls & C,
\end{align}
where the constant $C$ depends on $\|\rho_0\|_{L^\infty(0,T;H^m(\mathbb{T}))} $ and $\|\Om_0\|_{L^\infty(0,T;H^m(\mathbb{T}))} $ with $m>4$.

\noindent\textbf{Step IV}: Performing a similar scheme, we can get the boundedness of $\nm{\nabla_\om \nabla_x \p_t f_1}$ (and $\nm{\nabla_x \p_t f_1}$) with $m>7$, which yields that, together with the previous estimates,
\begin{align}
  \|f_1\|_{H^1_{M_0}} + \|\p_t f_1\|_{H^1_{M_0}} \le C.
\end{align}

Proving the lemma for the estimates in spaces $H^2_{M_0}$ (with $m>10$) and $H^3_{M_0}$ (with $m>13$) proceeds similarly as above, here we omit the details. \qed



\begin{remark}\label{Rema:auxiliary-high order}
  Actually, we can get the similar results for higher order derivatives of $f_1$,
  \begin{align}
    \sup_{t\in[0,T]} (\|f_1\|_{H^N_{M_0}}(t) + \|\p_t f_1\|_{H^N_{M_0}}(t)) \le C,
  \end{align}
where $N\ge 2$ and $C$ depends on $H^m(\mathbb{T})$ ($m>4N$) norm of $(\rho_0, \Om_0)$.
\end{remark}

\section{Energy Estimates for Remainder Equation} 
\label{sec:energy_for_remainder_eq}


The proof of Theorem \ref{Thm: main} relies crucially on the following a priori estimate for the remainder equation \eqref{Remainder-Equation-1}.
\begin{align}\label{Eq: Target} 
  \p_t (\ep f) + \om \cdot \nabla_x (\ep f) +\eta_0 \nabla_\om \cdot \left( P_{\om^\bot} D_{0} (\ep f) \right)
  + \mathcal{L}_0 f
  = h_1 + \ep f h_0,
\end{align}

We first introduce some instant energy functionals,
  \begin{align}
  \mathcal{E}= \mathcal{F}_0 +\mathcal{F}_1 +\mathcal{F}_2, \ \
  \mathcal{G}= \mathcal{G}_0 +\mathcal{G}_1 +\mathcal{G}_2, \ \
  \mathcal{H}= \mathcal{H}_0 +\mathcal{H}_1 +\mathcal{H}_2,
  \end{align}
where the terms $\mathcal{F}_i, \mathcal{G}_i, \mathcal{H}_i$ ($i=0,1,2$) are defined as follows,
\begin{align}
  &\mathcal{F}_0= \ep \nm{f}^2,\quad \mathcal{G}_0= \nm{\nabla_\om f}^2,\quad \mathcal{H}_0= \nm{h_1}^2,
  \\
  &\mathcal{F}_1=\ep^2 \nm{\nabla_x f}^2 + \ep \nm{\nabla_\om f}^2, \quad
  \mathcal{G}_1= \ep \nm{\nabla_\om \nabla_x f}^2 + \nm{\nabla_\om^2 f}^2, \\\no
  &\mathcal{H}_1= \ep \nm{\nabla_x h_1}^2 + \nm{\nabla_\om h_1}^2,
  \\
  &\mathcal{F}_2=\ep^3 \nm{\nabla_x^2 f}^2 + \ep^2 \nm{\nabla_\om \nabla_x f}^2 + \ep \nm{\nabla_\om^2 f}^2,\\\no
  &\mathcal{G}_2= \ep^2 \nm{\nabla_\om \nabla_x^2 f}^2 + \ep \nm{\nabla_\om^2 \nabla_x f}^2 + \nm{\nabla_\om^3 f}^2,\\\no
  &\mathcal{H}_2= \ep^2 \nm{\nabla_x^2 h_1}^2 + \ep \nm{\nabla_\om \nabla_x h_1}^2 + \nm{\nabla_\om^2 h_1}^2.
\end{align}

We shall give a uniform (in $\ep$) estimate for the instant energy $\mathcal{E}$, which plays the most important role in the process of Hilbert expansion.
\begin{lemma}\label{lemm: a priori}
  As for the above instant energy functional, there exists an $\ep_0>0$ such that, for any $\ep \in (0,\ep_0)$,
\begin{align}\label{eq: a priori}
\frac{\rm d}{{\rm d} t} \mathcal{E} + \mathcal{G}
\le \widetilde{C} (\mathcal{E} + \mathcal{H}),
\end{align}
where the constant $\widetilde{C}$ depends on the value of $\|\rho_0\|_{L^\infty(0,T; H^m(\mathbb{T}))}$ and $\|\Om_0\|_{L^\infty(0,T; H^m(\mathbb{T}))}$ (with $m>13$).
\end{lemma}


{\noindent\it Proof of Lemma \ref{lemm: a priori}}: we shall give three type energy estimates in what follows: $L^2_{M_0}$, $H^1_{M_0}$, $H^2_{M_0}$ estimates, corresponding to estimates for energy functionals $\mathcal{F}_0, \mathcal{F}_2, \mathcal{F}_2$, respectively, and then close them to finish the proof. We also give two remarks, one of which concerns the $H^N_{M_0}$ estimates (for large integer $N$) and the other discusses the role of $\mathcal{H}$.

\subsection{\rm Claim 1: ($L^2_{M_0}$ estimates).}
\label{sub:Claim 1}

Taking $L^2(M_0 \d \om \d x)$ inner product with $f$, we get from the equation \eqref{Eq: Target} that
%
%
%
%
%
%
\begin{align*}
\ep \la\la \p_t f, f \ra\ra_{_{M_0}}
= \frac{1}{2} \frac{\d}{\d t} \ep \|f\|^2_{L^2_{M_0}} - \frac{1}{2} \ep \iint |f|^2 \p_t M_0 \d \om \d x
\end{align*}
with \qquad $\displaystyle \left| \iint |f|^2 \p_t M_0 \d \om \d x \right| = \left| \iint |f|^2 \p_t \Om_0 \cdot \om M_0 \d \om \d x \right|
  \le \|\p_t \Om_0\|_{L^\infty_x} \|f\|^2_{M_0}.$

The second term on the left-hand side can be controlled by
\begin{align*}
  \la\la \om \cdot \nabla_x f, f \ra\ra_{M_0}
  =& \frac{1}{2} \iint \om \cdot \nabla_x (f^2 M_0) \d\om \d x
    - \frac{1}{2} \iint \om \cdot \nabla_x M_0 |f|^2 \d \om \d x \\
  =& 0 -\frac{1}{2} \iint \om \otimes \om: \nabla_x \Om_0 |f|^2 M_0 \d \om \d x,
\end{align*}
hencely $$\displaystyle \ep \la\la \om \cdot \nabla_x f, f \ra\ra_{M_0} \ls \ep \|\nabla_x \Om_0 \|_{L^\infty_x} \|f\|^2_{M_0}. $$

As for the third term, it follows that
\begin{align*}
  -\eta_0 \la\la \nabla_\om \cdot (P_{\om^\bot} D_0 f), f \ra\ra_{M_0} = & \eta_0 \la\la D_0 f, \nabla_\om f \ra\ra_{M_0}
 = \eta_0  \iint D_0  \frac{1}{2} \nabla_\om f^2 M_0 \d \om \d x \\
 = & \frac{1}{2} \eta_0 (n-1) \iint D_0 \om f^2 M_0 \d \om \d x - \frac{1}{2} \eta_0 \iint D_0 f^2 \nabla_\om M_0 \d \om \d x,
\end{align*}
combining with the facts that $\nabla_\om M_0= P_{\om^\bot} \Om_0 M_0$ and $|P_{\om^\bot} \Om_0| \le 1$, this gives
\begin{align*}
  |\eta_0 \ep \la\la \nabla_\om \cdot (P_{\om^\bot} D_0 f), f \ra\ra_{M_0} |
 \ls \ep \|D_0\|_{L^\infty_x} \|f\|^2_{M_0}.
\end{align*}

On the other hand, the nonnegativity of the linearized operator $\mathcal{L}_0 $ implies that
\begin{align*}
  \la\la \mathcal{L}_0 f, f \ra\ra_{M_0} = \|\nabla_\om f\|^2_{M_0}.
\end{align*}

The Poincar\'e inequality and the H\"older inequality gives the following estimate,
\begin{align*}
  \la\la h_1, f \ra\ra_{M_0} \le \|h_1\|_{M_0} \|f\|_{M_0} \ls \|h_1\|_{M_0} \|\nabla_\om f\|_{M_0}
 \le \frac{1}{2} C \|h_1\|^2_{M_0} + \frac{1}{2} \|\nabla_\om f\|^2_{M_0}.
\end{align*}

Notice that the last term on the right hand of equation \eqref{Eq: Target} can be bounded as follows,
\begin{align*}
  \ep \la\la h_0 f, f \ra\ra_{M_0} \le \|h_0\|_{M_0}\ \ep\|f\|^2_{M_0}
  \ls (\|\p_t \Om_0\|_{L^\infty_x} + \|\nabla_x \Om_0\|_{L^\infty_x} + \eta_0 \|D_0\|_{L^\infty_x}) \ep\|f\|^2_{M_0},
\end{align*}
then all these above estimates together gives that
\begin{align}
  \frac{1}{2} \frac{\d}{\d t} \ep\|f\|^2_{M_0} + \frac{1}{2} \|\nabla_\om f\|^2_{M_0}
 \ls (\|\p_t \Om_0\|_{L^\infty_x} + \|\nabla_x \Om_0\|_{L^\infty_x} + \eta_0 \|D_0\|_{L^\infty_x}) \ep\|f\|^2_{M_0}
     + \|h_1\|^2_{M_0}.
\end{align}

Denote
$$C_0 = \|\p_t \Om_0\|_{L^\infty_x} + \|\nabla_x \Om_0\|_{L^\infty_x} + \eta_0 \|D_0\|_{L^\infty_x},$$
then there exists a universal constant $C$ such that
\begin{align}\label{Esm: L^2 estimate}
  \frac{1}{2} \frac{\d}{\d t} \mathcal{F}_0 + \frac{1}{2} \mathcal{G}_0 \le C(C_0 \mathcal{F}_0 + \mathcal{H}_0).
\end{align}




\subsection{\rm Claim 2: ($H^1_{M_0} $ estimates)}
\label{sub:Claim 2}

\quad
\newline\textbf{Step I}: Taking $\m{L2} $ scalar product with $\mathcal{L}_0 f$ to control the quantity $\ep \norm*{M_0}{\nabla_\om f}^2$, we have the following estimates for the first two terms on the left-hand side,
\begin{align*}
  \ep \skp{\p_t f}{\mathcal{L}_0 f}
 =& \ep \skp{\p_t \nabla_\om f}{\nabla_\om f} \\
 =& \frac{1}{2} \frac{\d}{\d t} \ep \norm*{M_0}{\nabla_\om f}^2
   -\frac{1}{2}\ep  \iint \p_t M_0 |\nabla_\om f|^2 \d \om \d x \\
 =& \frac{1}{2} \frac{\d}{\d t} \ep \norm*{M_0}{\nabla_\om f}^2
   -\frac{1}{2}\ep  \iint \om \cdot \p_t \Om_0\ M_0 |\nabla_\om f|^2 \d \om \d x,
  \\
 \ep \skp{\om\cdot \nabla_x f}{\mathcal{L}_0 f}
 =& \ep \skp{\nabla_\om(\om\cdot \nabla_x f)}{\nabla_\om f}\\
 =& \ep \skp{\om\cdot \nabla_x f}{\nabla_\om f}
     + \ep \skp{\nabla_x f}{\nabla_\om f}\\
 =& \frac{1}{2}\ep \iint \om \cdot \nabla_x (|\nabla_\om f|^2 M_0) \d\om \d x
    - \frac{1}{2}\ep \iint \om \cdot \nabla_x M_0 |\nabla_\om f|^2 \d\om \d x\\
 &  + \ep \skp{\nabla_x f}{\nabla_\om f}\\
 =& 0 - \frac{1}{2}\ep \iint \om \otimes \om: \nabla_x \Om_0\ M_0 |\nabla_\om f|^2 \d\om \d x
    + \ep \norm*{M_0}{\nabla_x f} \norm*{M_0}{\nabla_\om f}\\
 \le & \ep \norm{ify}{\nabla_x \Om} \norm*{M_0}{\nabla_\om f}^2
      + C \ep^2 \norm*{M_0}{\nabla_x f}^2 + \frac{1}{4} \norm*{M_0}{\nabla^2_\om f}^2,
\end{align*}
where we have used the H\"older inequality and the Poincar\'e inequality on the sphere in the last line.

As for the third term, we write
\begin{align*}
  \ep \eta_0 \skp{\nabla_\om \cdot(P_{\om^\bot} D_0 f)}{\mathcal{L}_0 f}
=\ep \eta_0 \skp{-(n-1) D_0 \om f + D_0 \nabla_\om f}{\mathcal{L}_0 f}.
\end{align*}
It is an easy matter to calculate,
\begin{align*}
  \scalprod{L2om}{\om f}{\mathcal{L}_0 f} =& \scalprod{L2om}{\nabla_\om(\om f)}{\nabla_\om f}
 = \scalprod{L2om}{f + \om \nabla_\om f}{\nabla_\om f}\\
 =& \frac{1}{2} \int (\nabla_\om f^2) M_0 \d\om + \int \om |\nabla_\om f|^2 M_0 \d\om \\
 =& \frac{n-1}{2} \int \om f^2 M_0 \d\om - \frac{1}{2} \int \nabla_\om M_0  f^2 \d\om + \int \om |\nabla_\om f|^2 M_0 \d\om.
\end{align*}

By Lemma \ref{lemm: Liu-ARMA}, we get the fact
\begin{align}
  \scalprod{L2om}{\mathcal{L}_0 f}{\nabla_\om f} \le C \norm{L2om}{\nabla_\om f}^2.
\end{align}

Combining with the Poincar\'e inequality on the sphere, the two above inequalities yield that
\begin{align*}
  & | \ep \eta_0 \skp{\nabla_\om \cdot(P_{\om^\bot} D_0 f)}{\mathcal{L}_0 f}| \\
 \ls & \ep\eta_0 \norm{ify}{D_0} \norm*{M_0}{\nabla_\om f}^2
   + \ep\eta_0 \frac{(n-1)^2}{2} \norm{ify}{D_0} \norm*{M_0}{f}^2\\
   &+ \ep\eta_0 \frac{(n-1)}{2} \norm{ify}{D_0} \norm*{M_0}{f}^2
   + \ep\eta_0(n-1) \norm{ify}{D_0} \norm*{M_0}{\nabla_\om f}^2\\
 \ls & \ep \eta_0 \norm{ify}{D_0} \norm*{M_0}{\nabla_\om f}^2.
\end{align*}

On the other hand, we can control the terms on the right-hand side of equation \eqref{Eq: Target} as follows,
\begin{align*}
  \skp{\mathcal{L}_0 f}{\mathcal{L}_0 f}=& \nm{\nabla_\om^2 f}^2,
\\
  \skp{h_1}{\mathcal{L}_0 f}=& \skp{\nabla_\om h_1}{\nabla_\om f}
     \le \nm{\nabla_\om h_1} \nm{\nabla_\om f}
     \le C \nm{\nabla_\om h_1}^2 + \frac{1}{4} \nm{\nabla_\om^2 f}^2,
\\
  \ep \skp{h_0 f}{\mathcal{L}_0 f}=& \ep \skp{\nabla_\om f h_0}{\nabla_\om f} + \ep \skp{f \nabla_\om h_0}{\nabla_\om f}\\
    \le& \ep \nfyom{h_0} \nm{\nabla_\om f}^2 + \ep \nfyom{\nabla_\om h_0} \nm{f} \nm{\nabla_\om f},
\end{align*}
where we have used the Poincar\'e inequality again.

Accordingly, we have
\begin{align}\label{Esm:Deri-om}
   &\frac{1}{2}\frac{\d}{\d t} (\ep \nm{\nabla_\om f}^2) + \frac{1}{2} \nm{\nabla_\om^2 f}^2 \\\no
 \ls& \nm{\nabla_\om h_1}^2 + \ep^2 \nm{\nabla_x f}^2
    + (\nfy{\p_t \Om_0} + \nfy{\nabla_x \Om_0} + \eta_0 \nfy{D_0} )\ \ep \nm{\nabla_\om f}^2 \\\no
    &+ \ep \nfyom{h_0} \nm{\nabla_\om f}^2 + \ep \nfyom{\nabla_\om h_0} \nm{f} \nm{\nabla_\om f}.
\end{align}

\noindent\textbf{Step II}: Apply the operator $\nabla_x$ to the equation \eqref{Eq: Target}, and take $\m{L2} $ scalar product with $\ep \nabla_x f$, then we get
\begin{align*}
  \ep^2 \skp{\p_t \nabla_x f}{\nabla_x f}
  =& \frac{1}{2} \frac{\d}{\d x} \ep^2 \nm{\nabla_x f}^2 -\frac{1}{2} \ep^2 \iint \p_t M_0 |\nabla_x f|^2 \d\om \d x,
  \\
  \ep^2 \skp{\nabla_x(\om \cdot \nabla_x f)}{\nabla_x f}
  =& \frac{1}{2} \ep^2 \iint \om \cdot \nabla_x (|\nabla_x f|^2 M_0) \d\om \d x \\
    &- \frac{1}{2} \ep^2 \iint \om \cdot \nabla_x M_0 |\nabla_x f|^2  \d\om \d x \\
  =& -\frac{1}{2} \ep^2 \iint \om \otimes \om: \nabla_x \Om_0 M_0 |\nabla_x f|^2  \d\om \d x.
\end{align*}
The third term can be decomposed into two terms,
\begin{align*}
  &-\ep^2\eta_0 \skp{\nabla_x \nabla_\om \cdot (P_{\om^\bot} D_0 f)}{\nabla_x f} \\
  =& -\ep^2\eta_0 \skp{\nabla_\om \cdot (P_{\om^\bot} D_0 \nabla_x f)}{\nabla_x f}
     -\ep^2\eta_0 \skp{\nabla_\om \cdot (P_{\om^\bot} \nabla_x D_0 f)}{\nabla_x f} \\
  =& \ep^2\eta_0 (I_1 + I_2).
\end{align*}

Some straightforward calculations enable us to get
\begin{align*}
  I_1=& \la\! \scalprod{} {D_0 \nabla_x f}{\nabla_\om (\nabla_x f M_0)} \!\ra \\
     =& \la\! \scalprod{} {D_0 \nabla_x f}{\nabla_\om \nabla_x f M_0 + \nabla_x f \nabla_\om M_0} \!\ra \\
     =& \frac{1}{2} \iint D_0 (\nabla_\om |\nabla_x f|^2)\ M_0 \d\om \d x
       + \iint D_0 |\nabla_x f|^2 \nabla_\om M_0 \d\om \d x \\
     =&  \frac{1}{2} \iint D_0 (\nabla_\om |\nabla_x f|^2)\ M_0 \d\om \d x
       + \frac{1}{2} \iint D_0 |\nabla_x f|^2 M_0 P_{\om^\bot} \Om_0 \d\om \d x \\
     =& \frac{n-1}{2} \iint D_0 \om |\nabla_x f|^2 M_0 \d\om \d x
       + \frac{1}{2} \iint D_0 |\nabla_x f|^2 M_0 P_{\om^\bot} \Om_0 \d\om \d x \\
    \le & \frac{n}{2} \nfy{D_0} \nm{\nabla_x f}^2,
\end{align*}
and by the Poincar\'e inequality,
\begin{align*}
  I_2=& \skp{(n-1)\nabla_x D_0 \om f - \nabla_x D_0 \nabla_\om f}{\nabla_x f}\\
    \le & (n-1) \nfy{\nabla_x D_0} \nm{f} \nm{\nabla_x f}
        + \nfy{\nabla_x D_0} \nm{\nabla_\om f} \nm{\nabla_x f}\\
    \ls & n \nfy{\nabla_x D_0} \nm{\nabla_\om f} \nm{\nabla_x f},
\end{align*}
hence the third term can be bounded by
\begin{align*}
  &|\ep^2\eta_0 \skp{\nabla_x \nabla_\om \cdot (P_{\om^\bot} D_0 f)}{\nabla_x f} |
 = \ep^2\eta_0 |I_1 + I_2|\\
 \ls& \ep^2\eta_0 \nfy{D_0} \nm{\nabla_x f}^2
    + \ep^2\eta_0 \nfy{\nabla_x D_0} \nm{\nabla_\om f} \nm{\nabla_x f}.
\end{align*}

Next we turn to consider the terms on the right-hand side of equation \eqref{Eq: Target}. Recalling the nonnegativity property of the linearized operator $\mathcal{L}_0=-(\Delta_\om + \Om_0 \cdot \nabla_\om)$, we have
\begin{align*}
  \skp{\nabla_x \mathcal{L}_0 f}{\ep \nabla_x f}
 =&\ep \skp{\mathcal{L}_0 \nabla_x f}{\nabla_x f}
   -\ep \skp{\nabla_x \Om_0 \nabla_\om f}{\nabla_x f}\\
 \ge &\ep \nm{\nabla_\om \nabla_x f}^2
  - \ep \nfy{\nabla_x \Om_0} \nm{\nabla_\om f} \nm{\nabla_x f} \\
 \ge & \frac{3}{4}\ep \nm{\nabla_\om \nabla_x f}^2
  - C \ep \nfy{\nabla_x \Om_0}^2 \nm{\nabla_\om f}^2,
\end{align*}
where we have used the Poincar\'e inequality on the sphere $\nm{\nabla_x f}^2 \ls \nm{\nabla_\om \nabla_x f}^2 $.

Note the Poincar\'e inequality also yield that
\begin{align*}
  \skp{\nabla_x h_1}{\ep \nabla_x f}
 \le& \ep \nm{\nabla_x h_1} \nm{\nabla_x f}
 \le C \ep \nm{\nabla_x h_1}^2 + \frac{1}{4}\ep \nm{\nabla_\om \nabla_x f}^2,
 \\
  \ep^2 \skp{\nabla_x(h_0 f)}{\nabla_x f}
  =& \ep^2 \skp{\nabla_x f h_0}{\nabla_x f} + \ep^2 \skp{f \nabla_x h_0}{\nabla_x f} \\
  \le& \ep^2 \nfyom{h_0} \nm{\nabla_x f}^2 + \ep^2 \nfyom{\nabla_x h_0} \nm{f} \nm{\nabla_x f},
\end{align*}
which ensure that, together with the previous estimates,
\begin{align}\label{Esm:Deri-x}
  & \frac{1}{2}\frac{\d}{\d t} (\ep^2 \nm{\nabla_x f}^2) + \frac{1}{2} \ep \nm{\nabla_\om \nabla_x f}^2 \\\no
 \ls& \ep \nfy{\nabla_x \Om_0}^2 \nm{\nabla_\om f}^2
    +\ep \nm{\nabla_x h_1}^2 \\\no
    & + (\nfy{\p_t \Om_0} + \nfy{\nabla_x \Om_0} + \eta_0 \nfy{D_0})\ \ep^2 \nm{\nabla_x f}^2 \\\no
    & + \ep^2 \eta_0 \nfy{\nabla_x D_0} \nm{\nabla_\om f} \nm{\nabla_x f} \\\no
    & + (\ep^2 \nfyom{h_0} \nm{\nabla_x f}^2
      + \ep^2 \nfyom{\nabla_x h_0} \nm{f} \nm{\nabla_x f} ).
\end{align}

\noindent\textbf{Step III}: Combining the two inequalities \eqref{Esm:Deri-om} and \eqref{Esm:Deri-x} gives that
\begin{align}
  & \frac{1}{2}\frac{\d}{\d t} \lt(\ep^2 \nm{\nabla_x f}^2
  + \ep \nm{\nabla_\om f}^2 \rt) + \frac{1}{2} \lt(\ep \nm{\nabla_\om \nabla_x f}^2 + \nm{\nabla_\om^2 f}^2 \rt) \\\no
 \ls & \lt(\ep^2 \nm{\nabla_x f}^2 + \nfy{\nabla_x \Om_0}^2\ \ep \nm{\nabla_\om f}^2 \rt)
  + \lt(\ep \nm{\nabla_x h_1}^2 + \nm{\nabla_\om h_1}^2 \rt) \\\no
  &+ \lt(\nfy{\p_t \Om_0} + \nfy{\nabla_x \Om_0} + \eta_0 \nfy{D_0}\rt) (\ep^2 \nm{\nabla_x f}^2 + \ep \nm{\nabla_\om f}^2) \\\no
  & + \eta_0 \nfy{\nabla_x D_0}\ \ep^2 \nm{\nabla_\om f} \nm{\nabla_x f} \\\no
    & + \Big\{\nfyom{h_0}\ (\ep^2 \nm{\nabla_x f}^2 + \ep \nm{\nabla_\om f}^2)
      + \nfyom{\nabla_x h_0} \ep^2 \nm{f} \nm{\nabla_x f} \\\no
    & + \nfyom{\nabla_\om h_0}\ \ep \nm{f} \nm{\nabla_\om f} \Big\}.
\end{align}

Denote $$C_1 =\sup\lt\{C_0^2,\ \nfyom{\nabla_x h_0},\ \nfyom{\nabla_\om h_0},\ \eta_0 \nfy{\nabla_x D_0}\rt\},$$
with $C_0 = \nfy{\p_t \Om_0} + \nfy{\nabla_x \Om_0} + \eta_0 \nfy{D_0}$ defined as before, then we can conclude that, up to a constant $C>0$,
\begin{align}\label{Esm: H^1 estimate}
  \frac{1}{2} \frac{\d}{\d t} \mathcal{F}_1 + \frac{1}{2} \mathcal{G}_1
 \le C C_1 (1+ \ep^\frac{1}{2}) (\mathcal{F}_1 + \mathcal{F}_0) + C \mathcal{H}_1,
\end{align}
where it should be pointed out that
\begin{align*}
  \nfyom{h_0} \le& C_0,
\\
  \ep^2 \nm{\nabla_\om f} \nm{\nabla_x f}
 =& \ep^\frac{1}{2} (\ep^\frac{1}{2}\nm{\nabla_\om f}) (\ep \nm{\nabla_x f})
   \le \ep^\frac{1}{2}\mathcal{F}_1^\frac{1}{2} \mathcal{F}_1^\frac{1}{2}.
\end{align*}



\subsection{\rm Claim 3: ($H^2_{M_0} $ estimates)}
\label{sub:Claim 3}

\quad
\newline\textbf{Step I}: Apply the operator $\nabla_x^2 $ to the equation \eqref{Eq: Target}, and take $\m{L2} $ scalar product with the quantity $\ep^2 \nabla_x^2 f$, then we can get immediately,
\begin{align*}
  \ep^3 \skp{\p_t \nabla_x^2 f}{\nabla_x^2 f}
 =& \frac{1}{2} \frac{\d }{\d t} \ep^3 \nm{\nabla_x^2 f}^2
   - \frac{1}{2} \ep^3 \iint \p_t M_0\ |\nabla_x f|^2 \d\om \d x,
  \\
  \ep^3 \skp{\nabla_x^2(\om \cdot \nabla_x f)}{\nabla_x^2 f}
 =& \frac{1}{2} \ep^3 \iint \om \cdot \nabla_x (|\nabla_x^2 f|^2) M_0 \d\om \d x \\
 =& 0 -\frac{1}{2} \ep^3 \iint \om \cdot \nabla_x M_0\ |\nabla_x^2 f|^2 \d\om \d x.
\end{align*}

Similar as before, we split the third term into two terms,
\begin{align*}
&-\ep^3 \eta_0 \skp{\nabla_x^2 \nabla_\om \cdot (P_{\om^\bot} D_0 f)}{\nabla_x^2 f} \\
=&-\ep^3 \eta_0 \skp{\nabla_\om \cdot [P_{\om^\bot} (\nabla_x^2 D_0 f + 2 \nabla_x D_0 \nabla_x f)]}{\nabla_x^2 f} \\
&-\ep^3 \eta_0 \skp{\nabla_\om \cdot (P_{\om^\bot} D_0 \nabla_x^2 f)}{\nabla_x^2 f} \\
=& \ep^3 \eta_0 (II_1 + II_2).
\end{align*}

It follows by calculation that,
\begin{align*}
  II_1=& \skp{(n-1)\nabla_x^2 D_0 \om f - \nabla_x^2 D_0 \nabla_\om f}{\nabla_x^2 f} \\
  &+ 2 \skp{(n-1)\nabla_x D_0 \om \nabla_x f - \nabla_x D_0 \nabla_\om \nabla_x f}{\nabla_x^2 f} \\
  \le & \nfy{\nabla_x^2 D_0} \lt( (n-1) \nm{f} + \nm{\nabla_\om f} \rt)  \nm{\nabla_x^2 f} \\
    & + 2 \nfy{\nabla_x D_0} \lt( (n-1) \nm{\nabla_x f} + \nm{\nabla_\om \nabla_x f} \rt)  \nm{\nabla_x^2 f} \\
  \ls & \nfy{\nabla_x^2 D_0} \nm{\nabla_\om f} \nm{\nabla_x^2 f} + \nfy{\nabla_x D_0} \nm{\nabla_\om \nabla_x f} \nm{\nabla_x^2 f},
\\
  II_2=& \skp{D_0 \nabla_x^2 f}{\nabla_\om \nabla_x^2 f}
  =\frac{1}{2} \iint D_0 \nabla_\om (|\nabla_x^2 f|^2)\ M_0 \d\om \d x \\
  =& \frac{n-1}{2} \iint D_0 \om\ |\nabla_x^2 f|^2\ M_0 \d\om \d x
     -\frac{1}{2} \iint D_0 \nabla_\om M_0\ |\nabla_x^2 f|^2 \d\om \d x \\
  \ls & \nfy{D_0} \nm{\nabla_x^2 f}^2,
\end{align*}
which give the estimate
\begin{align*}
&\lt|\ep^3 \eta_0 \skp{\nabla_x^2 \nabla_\om \cdot (P_{\om^\bot} D_0 f)}{\nabla_x^2 f} \rt| \\
\ls & \ep^3 \eta_0 \nfy{D_0} \nm{\nabla_x^2 f}^2
    + \ep^3 \eta_0 \nfy{\nabla_x D_0} \nm{\nabla_\om \nabla_x f} \nm{\nabla_x^2 f}\\
  & + \ep^3 \eta_0 \nfy{\nabla_x^2 D_0} \nm{\nabla_\om f} \nm{\nabla_x^2 f}.
\end{align*}

As for the right-hand side of the equation \eqref{Eq: Target}, we collect their estimates, by the nonnegativity property of the linearized operator $\mathcal{L}_0$ and the Poincar\'e inequality on the sphere $\nm{\nabla_x^2 f} \ls \nm{\nabla_\om \nabla_x^2 f}$, as follows,
\begin{align*}
  &\skp{\nabla_x^2 \mathcal{L}_0 f}{\ep^2 \nabla_x^2 f} \\
  =& \ep^2 \skp {\mathcal{L}_0 \nabla_x^2 f}{\nabla_x^2 f} - \ep^2 \skp{\nabla_x^2 \Om_0 \nabla_\om f + 2 \nabla_x \Om_0 \nabla_\om \nabla_x f}{\nabla_x^2 f} \\
  \ge & \ep^2 \nm{\nabla_\om \nabla_x^2 f}^2
       - \Big[ \nfy{\nabla_x^2 \Om_0} \nm{\nabla_\om f} \nm{\nabla_x^2 f} \\
     & + 2 \nfy{\nabla_x \Om_0} \nm{\nabla_\om \nabla_x f} \nm{\nabla_x^2 f} \Big] \\
  \ge & \frac{3}{4}\ep^2 \nm{\nabla_\om \nabla_x^2 f}^2
     - C \ep^2 \lt( \nfy{\nabla_x^2 \Om_0}^2 \nm{\nabla_\om f}^2 + \nfy{\nabla_x \Om_0}^2 \nm{\nabla_\om \nabla_x f}^2 \rt),
\\
   &\ep^2 \skp{\nabla_x^2 h_1}{\nabla_x^2 f} \\
  \le &\ep^2 \nm{\nabla_x^2 h_1} \nm{\nabla_x^2 f}
  \ls \ep^2 \nm{\nabla_x^2 h_1} \nm{\nabla_\om \nabla_x^2 f} \\
  \le & C \ep^2 \nm{\nabla_x^2 h_1}^2 + \frac{1}{4}\ep^2 \nm{\nabla_\om \nabla_x^2 f}^2,
\\
  & \ep^3 \skp{\nabla_x^2 (h_0 f)}{\nabla_x^2 f} \\
  \ls & \ep^3 \nfyom{h_0} \nm{\nabla_x^2 f}^2 + \ep^3 \nfyom{\nabla_x h_0} \nm{\nabla_x f} \nm{\nabla_x^2 f} \\
      & + \ep^3 \nfyom{\nabla_x^2 h_0} \nm{f} \nm{\nabla_x^2 f}.
\end{align*}

Then it follows from the above estimates that,
\begin{align}\label{Esm:Deri-x2}
  &\frac{1}{2} \frac{\d}{\d t} \ep^3 \nm{\nabla_x^2 f}^2 + \frac{1}{2} \ep^2 \nm{\nabla_\om \nabla_x^2 f}^2 \\\no
 \ls& \ep^2 ( \nfy{\nabla_x^2 \Om_0}^2 \nm{\nabla_\om f}^2 + \nfy{\nabla_x \Om_0}^2 \nm{\nabla_\om \nabla_x f}^2 ) + \ep^2 \nm{\nabla_x^2 h_1}^2 \\\no
 & + (\nfy{\p_t \Om_0} + \nfy{\nabla_x \Om_0} + \eta_0 \nfy{D_0} )\ \ep^3 \nm{\nabla_x^2 f}^2 \\\no
 & + \ep^3 \eta_0 \nfy{\nabla_x D_0} \nm{\nabla_\om \nabla_x f} \nm{\nabla_x^2 f}  
  + \ep^3 \eta_0 \nfy{\nabla_x^2 D_0} \nm{\nabla_\om f} \nm{\nabla_x^2 f}  \\\no
 & + \Big(\nfyom{h_0} \ep^3 \nm{\nabla_x^2 f}^2 + \nfyom{\nabla_x h_0} \ep^3 \nm{\nabla_x f} \nm{\nabla_x^2 f} \\\no
 & + \nfyom{\nabla_x^2 h_0} \ep^3 \nm{f} \nm{\nabla_x^2 f}\Big).
\end{align}

\noindent\textbf{Step II}: To control the quantity $\ep^2 \nm{\nabla_\om \nabla_x f}^2$, we apply the operator $\nabla_x$ to the equation \eqref{Eq: Target}, and take $\m{L2}$ scalar product with $\ep \mathcal{L}_0 \nabla_x f$, then it follows
\begin{align*}
  & \ep^2 \skp{\p_t \nabla_x f}{\mathcal{L}_0 \nabla_x f} \\
  =& \ep^2 \skp{\p_t \nabla_\om \nabla_x f}{\nabla_\om \nabla_x f} \\
  =& \frac{1}{2} \frac{\d}{\d t} \ep^2 \nm{\nabla_\om \nabla_x f}^2 - \frac{1}{2} \ep^2 \iint \p_t M_0 |\nabla_\om \nabla_x f|^2 \d\om \d x,
\\[2pt]
  & \ep^2 \skp{\om\cdot \nabla_x^2 f}{\mathcal{L}_0 \nabla_x f} \\
  =& \ep^2 \skp{\nabla_\om (\om \cdot \nabla_x^2 f)}{\nabla_\om \nabla_x f} \\
  =& \frac{1}{2} \ep^2 \iint \om \cdot \nabla_x (|\nabla_\om \nabla_x f|^2 M_0) \d\om \d x
    -\frac{1}{2} \ep^2 \iint \om \cdot \nabla_x M_0 |\nabla_\om \nabla_x f|^2 \d\om \d x \\
    & + \ep^2 \skp{\nabla_x^2 f}{\nabla_\om \nabla_x f} \\
  \le & 0 + \frac{1}{2} \ep^2 \nfy{\nabla_x \Om_0} \nm{\nabla_\om \nabla_x f}^2
      + C \ep^3 \nm{\nabla_x^2 f}^2 + \frac{1}{6} \ep \nm{\nabla_\om^2 \nabla_x f}^2,
\end{align*}
where we have used the Poincar\'e inequality on the sphere $\nm{\nabla_\om \nabla_x f} \ls \nm{\nabla_\om^2 \nabla_x f}$.

The third term can be estimated by a similar decomposition as before,
\begin{align*}
    & -\ep^2 \eta_0 \skp{\nabla_x \nabla_\om \cdot (P_{\om^\bot} D_0 f)}{\mathcal{L}_0 \nabla_x f} \\
  = & -\ep^2 \eta_0 \skp{\nabla_\om \cdot (P_{\om^\bot} D_0 \nabla_x f + P_{\om^\bot} \nabla_x D_0 f) }{\mathcal{L}_0 \nabla_x f} \\
  = & \ep^2 \eta_0 \skp{(n-1)D_0 \om \nabla_x f - D_0 \nabla_\om \nabla_x f}{\mathcal{L}_0 \nabla_x f} \\
     &+ \ep^2 \eta_0 \skp{(n-1) \nabla_x D_0 \om f - \nabla_x D_0 \nabla_\om f}{\mathcal{L}_0 \nabla_x f} \\
  = & \ep^2 \eta_0 (II_3 + II_4).
\end{align*}

Combining Lemma \ref{lemm: Liu-ARMA} and the Poincar\'e inequality on the sphere ensures that
\begin{align*}
  II_3 =& (n-1)\skp{D_0 \nabla_\om(\om \cdot \nabla_x f)}{\nabla_\om \nabla_x f} - \skp{D_0 \nabla_\om \nabla_x f}{\mathcal{L}_0 \nabla_x f} \\
  \le & (n-1) \nfy{D_0} (\nm{\nabla_x f} + \nm{\nabla_\om \nabla_x f}) \nm{\nabla_\om \nabla_x f} \\
    & + \iint D_0 \mathcal{L}_0 \nabla_x f \nabla_\om \nabla_x f M_0 \d\om \d x\\
  \ls & (n-1) \nfy{D_0} \nm{\nabla_\om \nabla_x f}^2 + \nfy{D_0} \nm{\nabla_\om \nabla_x f}^2\\
  \le & n \nfy{D_0} \nm{\nabla_\om \nabla_x f}^2,
\end{align*}
and
\begin{align*}
  II_4 =& (n-1)\skp{\nabla_x D_0 \nabla_\om (\om f)}{\nabla_\om \nabla_x f} - \skp{\nabla_x D_0 \nabla_\om^2 f}{\nabla_\om \nabla_x f} \\
  \le & (n-1) \nfy{\nabla_x D_0} (\nm{f} + \nm{\nabla_\om f}) \nm{\nabla_\om \nabla_x f} \\
   &+ \nfy{\nabla_x D_0} \nm{\nabla_\om^2 f} \nm{\nabla_\om \nabla_x f} \\
  \ls & n \nfy{\nabla_x D_0} \nm{\nabla_\om^2 f} \nm{\nabla_\om \nabla_x f},
\end{align*}
therefore we get
\begin{align*}
  &\lt|\ep^2 \eta_0 \skp{\nabla_x \nabla_\om \cdot (P_{\om^\bot} D_0 f)}{\mathcal{L}_0 \nabla_x f}\rt|\\
  \ls& \ep^2 \eta_0 \nfy{D_0} \nm{\nabla_\om \nabla_x f}^2 + \ep^2 \eta_0 \nfy{\nabla_x D_0} \nm{\nabla_\om^2 f} \nm{\nabla_\om \nabla_x f}.
\end{align*}

Similar discussions as before enable us to deal with the terms on the right-hand side, as follows,
\begin{align*}
  \ep \skp{\nabla_x \mathcal{L}_0 f}{\mathcal{L}_0 \nabla_x f}
  = & \ep \skp{\mathcal{L}_0 \nabla_x f}{\mathcal{L}_0 \nabla_x f}
     - \ep \skp{\nabla_x \Om_0 \nabla_\om f}{\mathcal{L}_0 \nabla_x f} \\
  \ge & \ep \nm{\nabla_\om^2 \nabla_x f}^2
      - \ep \nfy{\nabla_x \Om_0} \nm{\nabla_\om^2 f} \nm{\nabla_\om \nabla_x f} \\
  \ge & \frac{5}{6} \ep \nm{\nabla_\om^2 \nabla_x f}^2
      - C \ep \nfy{\nabla_x \Om_0}^2 \nm{\nabla_\om^2 f}^2,
\\
  \ep \skp{\nabla_x h_1}{\mathcal{L}_0 \nabla_x f} =& \ep \skp{\nabla_\om \nabla_x h_1}{\nabla_\om \nabla_x f}
  \le C \ep \nm{\nabla_\om \nabla_x h_1}^2 + \frac{1}{6} \ep \nm{\nabla_\om^2 \nabla_x f}^2,
\\
  \ep^2 \skp{\nabla_x(h_0 f)} {\mathcal{L}_0 \nabla_x f} =& \ep^2 \skp{\nabla_\om \nabla_x(h_0 f)}{\nabla_\om \nabla_x f} \\
  \ls & \ep^2 \sum_{\substack{|\alpha_1|+|\alpha_2|=1\\ |\beta_1|+|\beta_2|=1}} \nfyom{\nabla_\om^{\alpha_1} \nabla_x^{\beta_1} h_0}
     \nm{\nabla_\om^{\alpha_2} \nabla_x^{\beta_2} f} \nm{\nabla_\om \nabla_x f}.
\end{align*}

All the above estimates together give that
\begin{align}\label{Esm:Deri-om-x}
  & \frac{1}{2} \frac{\d}{\d t} \ep^2 \nm{\nabla_\om \nabla_x f}^2 + \frac{1}{2} \ep \nm{\nabla_\om^2 \nabla_x f}^2 \\\no
 \ls & \ep \nfy{\nabla_x \Om_0}^2 \nm{\nabla_\om^2 f}^2 + \ep \nm{\nabla_\om \nabla_x h_1}^2 + \ep^3 \nm{\nabla_x^2 f}^2 \\\no
 & + (\nfy{\p_t \Om_0} + \nfy{\nabla_x \Om_0} + \eta_0 \nfy{D_0} )\ \ep^2 \nm{\nabla_\om \nabla_x f}^2 \\\no
 & + \ep^2 \eta_0 \nfy{\nabla_x D_0} \nm{\nabla_\om^2 f} \nm{\nabla_\om \nabla_x f} \\\no
 & + \ep^2 \sum_{\substack{|\alpha_1|+|\alpha_2|=1\\ |\beta_1|+|\beta_2|=1}}
     \nfyom{\nabla_\om^{\alpha_1} \nabla_x^{\beta_1} h_0}
     \nm{\nabla_\om^{\alpha_2} \nabla_x^{\beta_2} f} \nm{\nabla_\om \nabla_x f}.
\end{align}

\noindent\textbf{Step III}: By taking $\m{L2}$ scalar product with $\mathcal{L}_0^2 f$, we will get from the equation \eqref{Eq: Target} the control for $\ep \nm{\nabla_\om^2 f}^2$ in a similar process as before. Here we omit the details and only list the estimates,
\begin{align*}
  \ep \skp{\p_t f}{\mathcal{L}_0^2 f}
 =\la\!\la & \p_t \nabla_\om^2 f, \nabla_\om^2 f \ra\!\ra_{M_0}  \\
 &\hspace*{-2em}
 =\frac{1}{2} \frac{\d}{\d t} \ep \nm{\nabla_\om^2 f}^2 - \frac{1}{2} \ep \iint \p_t M_0 |\nabla_\om^2 f|^2 \d\om \d x,
\\
 \ep \skp{\om \cdot \nabla_x f}{\mathcal{L}_0^2 f}
 =& \ep \skp{\om \cdot \nabla_x \nabla_\om^2  f}{\nabla_\om^2 f}
    + 2 \ep \skp{\nabla_\om \nabla_x f}{\nabla_\om^2 f} \\
 =& \frac{1}{2} \ep \iint \om \cdot \nabla_x (|\nabla_\om^2 f|^2 M_0) \d\om \d x
   - \frac{1}{2} \ep \iint \om \cdot \nabla_x M_0 |\nabla_\om^2 f|^2 \d\om \d x \\
  & + 2 \ep \skp{\nabla_\om \nabla_x f}{\nabla_\om^2 f} \\
 \le & \ep \nfy{\nabla_x \Om_0} \nm{\nabla_\om^2 f}^2 + C \ep^2 \nm{\nabla_\om \nabla_x f}^2 + \frac{1}{4} \nm{\nabla_\om^3 f}^2,
\\
  |\ep \eta_0 \la\!\la \nabla_\om \cdot (P_{\om^\bot} D_0 f), & \mathcal{L}_0^2 f \ra\!\ra_{M_0}|
 \le \ep \eta_0 \lt|\skp{(n-1) D_0 \om f - D_0 \nabla_\om f}{\mathcal{L}_0^2 f} \rt| \\
 \ls & \ep \eta_0 \nfy{D_0} \nm{\nabla_\om^2 f}^2, \ (\textrm{Using Lemma \ref{lemm: Liu-ARMA}})
\\
 \skp{\mathcal{L}_0 f}{\mathcal{L}_0^2 f} =& \nm{\nabla_\om^3 f}^2,
\\
 \skp{h_1}{\mathcal{L}_0^2 f} \le & \nm{\nabla_\om^2 h_1} \nm{\nabla_\om^2 f}
 \le C \nm{\nabla_\om^2 h_1}^2 + \frac{1}{4} \nm{\nabla_\om^3 f}^2,
\\
 \ep \skp{h_0 f}{\mathcal{L}_0^2 f}
 = & \ep \skp{\nabla_\om^2 (h_0 f)}{\nabla_\om^2 f} \\
 \ls & \ep \sum_{|\alpha_1| + |\alpha_2| =2} \nfyom{\nabla_\om^{\alpha_1} h_0} \nm{\nabla_\om^{\alpha_2} f} \nm{\nabla_\om^2 f} \\
 \ls & \ep \sum_{0\le |\alpha_1| \le 2} \nfyom{\nabla_\om^{\alpha_1} h_0} \nm{\nabla_\om^2 f}^2.
\end{align*}

Combining these above estimates implies that
\begin{align}\label{Esm:Deri-om2}
   &\frac{1}{2} \frac{\d}{\d t} \ep \nm{\nabla_\om^2 f}^2 + \frac{1}{2} \nm{\nabla_\om^3 f}^2 \\\no
 \ls& \nm{\nabla_\om^2 h_1}^2 + \ep^2 \nm{\nabla_\om \nabla_x f}^2 + \sum_{0\le |\alpha_1| \le 2} \nfyom{\nabla_\om^{\alpha_1} h_0} \ep \nm{\nabla_\om^2 f}^2 \\\no
 & + (\nfy{\p_t \Om_0} + \nfy{\nabla_x \Om_0} + \eta_0 \nfy{D_0})\ \ep \nm{\nabla_\om^2 f}^2.
\end{align}

\noindent\textbf{Step IV}: Summing up the inequalities \eqref{Esm:Deri-x2}, \eqref{Esm:Deri-om-x} and \eqref{Esm:Deri-om2}, we get
\begin{align}
  & \frac{1}{2} \frac{\d}{\d t} \lt(\ep^3 \nm{\nabla_x^2 f}^2 + \ep^2 \nm{\nabla_\om \nabla_x f}^2 + \ep \nm{\nabla_\om^2 f}^2 \rt) \\\no
   & + \frac{1}{2} \lt(\ep^2 \nm{\nabla_\om \nabla_x^2 f}^2 + \ep \nm{\nabla_\om^2 \nabla_x f}^2 + \nm{\nabla_\om^3 f}^2 \rt) \\\no
 \ls & \ep^2 ( \nfy{\nabla_x^2 \Om_0}^2 \nm{\nabla_\om f}^2 + \nfy{\nabla_x \Om_0}^2 \nm{\nabla_\om \nabla_x f}^2) + \ep \nfy{\nabla_x \Om_0}^2 \nm{\nabla_\om^2 f}^2 \\\no
   & + \lt(\ep^2 \nm{\nabla_x^2 h_1}^2 + \ep \nm{\nabla_\om \nabla_x h_1}^2 + \nm{\nabla_\om^2 h_1}^2 \rt) 
    + \ep^3 \nm{\nabla_x^2 f}^2 + \ep^2 \nm{\nabla_\om \nabla_x f}^2 \\\no
   & + (\nfy{\p_t \Om_0} + \nfy{\nabla_x \Om_0} + \eta_0 \nfy{D_0})
       (\ep^3 \nm{\nabla_x^2 f}^2 + \ep^2 \nm{\nabla_\om \nabla_x f}^2 + \ep \nm{\nabla_\om^2 f}^2 ) \\\no
   & + \eta_0 \Big[ \ep^3 \nfy{\nabla_x D_0} \nm{\nabla_\om \nabla_x f} \nm{\nabla_x^2 f} 
    + \ep^3 \nfy{\nabla_x^2 D_0} \nm{\nabla_\om f} \nm{\nabla_x^2 f} \\\no
  &\quad + \ep^2 \nfy{\nabla_x D_0} \nm{\nabla_\om^2 f} \nm{\nabla_\om \nabla_x f} \Big] \\\no
  & + \ep^3 \sum_{|\beta_1|+|\beta_2|=2} \nfyom{\nabla_x^{\beta_1} h_0} \nm{\nabla_x^{\beta_2} f} \nm{\nabla_x^2 f} \\\no
  & + \ep^2 \sum_{\substack{|\alpha_1|+|\alpha_2|=1\\ |\beta_1|+|\beta_2|=1}}
     \nfyom{\nabla_\om^{\alpha_1} \nabla_x^{\beta_1} h_0}
     \nm{\nabla_\om^{\alpha_2} \nabla_x^{\beta_2} f} \nm{\nabla_\om \nabla_x f}\\\no
  & + \ep \sum_{0\le |\alpha_1| \le 2} \nfyom{\nabla_\om^{\alpha_1} h_0} \nm{\nabla_\om^2 f}^2.
\end{align}

Denote $$C_2 =\sup\Big\{C_0^2, \ \nfy{\nabla_x^2 \Om_0}^2,\ \eta_0 \nfy{\nabla_x^2 D_0}^2,\
 \sum_{\mathclap{|\alpha|+|\beta|\le 2}} \nfyom{\nabla_\om^{\alpha} \nabla_x^{\beta} h_0}^2 \Big\},$$
then the above equation can be rewritten as
\begin{align}\label{Esm: H^2 estimate}
  \frac{1}{2} \frac{\d}{\d t} \mathcal{F}_2 + \frac{1}{2} \mathcal{G}_2
\le& C C_2 (1+\ep + \eta_0 \ep^\frac{1}{2} + \eta_0 \ep) \mathcal{F}_2
     + C \mathcal{H}_2 \\\no
   & + C C_2 (1+\ep^\frac{1}{2}) \mathcal{F}_2^\frac{1}{2} \mathcal{F}_1^\frac{1}{2}
     + C C_2 (1+\ep^\frac{1}{2} + \ep) \mathcal{F}_2^\frac{1}{2} \mathcal{F}_0^\frac{1}{2} \\\no
\le& C C_2 (1+\ep + \eta_0 \ep^\frac{1}{2} + \eta_0 \ep) \mathcal{F}_2
     + C C_2 (1+\ep^\frac{1}{2}) \mathcal{F}_1 \\\no
    &+ C C_2 (1+ \ep) \mathcal{F}_0
     + C \mathcal{H}_2 \\\no
\le& 2CC_2 (\mathcal{F}_2 + \mathcal{F}_1 + \mathcal{F}_0) + C \mathcal{H}_2.
\end{align}



\subsection{Closing the Estimates} 
\label{sub:closing_the_estimates}


Noticing the definition of $\mathcal{E}=\mathcal{F}_0 + \mathcal{F}_1 + \mathcal{F}_2$, we collect the inequalities \eqref{Esm: L^2 estimate}, \eqref{Esm: H^1 estimate}, and \eqref{Esm: H^2 estimate} to get that
\begin{align}
\frac{\d}{\d t} \mathcal{E} + \mathcal{G}
\le \widetilde{C} (\mathcal{E} + \mathcal{H}),
\end{align}
holds for any $\ep \in (0,\ep_0)$ and $t\in[0,T]$. Obviously, the value of $\widetilde{C}$ depends upon the value of $\sup_{t\in[0,T]} \{C_0, C_1, C_2\}$, which is depending eventually upon the value of $\|\rho_0\|_{L^\infty(0,T; H^m(\mathbb{T}))}$ and $\|\Om_0\|_{L^\infty(0,T; H^m(\mathbb{T}))}$ (with $m>13$). This completes the whole proof of Lemma \ref{lemm: a priori}. \qed


\begin{remark}\label{Rema:a priori-high order}
Supposed that the initial data $(\rho_0^{in}, \Om_0^{in})$ are sufficiently smooth, then we can consider high order energy estimates for the remainder equation \eqref{Eq: Target}. Here we only give the result about the $H^N_{x,\om}$ energy estimates because the proof is similar as before. Specifically speaking,
\begin{align}
  & \frac{1}{2} \frac{\rm d}{{\rm d} t} \sum_{\substack{|\alpha|+|\beta|=N\\ 0\le |\beta|\le N}} \ep^{|\beta|+1} \nm{\nabla_\om^\alpha \nabla_x^\beta f}^2  + \frac{1}{2} \sum_{\substack{|\alpha|+|\beta|=N\\ 0\le |\beta|\le N}} \ep^{|\beta|} \nm{\nabla_\om^{\alpha+1} \nabla_x^\beta f}^2 \\\no
 \ls & \sum_{\substack{|\alpha|+|\beta|=N\\ 1\le |\beta|\le N}}
       \sum_{\substack{\beta_1+\beta_2=\beta\\ |\beta_1|\ge 1}} C_{\beta}^{\beta_1}
        \nfy{\nabla_x^{\beta_1} \Om_0}^2 \nm{\nabla_\om^{\alpha+1} \nabla_x^{\beta_2} f}^2 \ep^{|\beta|}\\\no
   & + \sum_{\substack{|\alpha|+|\beta|=N\\ 0\le |\beta|\le N}} \ep^{|\beta|} \nm{\nabla_\om^{\alpha} \nabla_x^\beta h_1}^2
     + \sum_{\substack{|\alpha|+|\beta|=N\\ 0\le |\beta|\le N-1}} \ep^{|\beta|+2} \nm{\nabla_\om^{\alpha-1} \nabla_x^{\beta+1} f}^2 \\\no
   & + (\nfy{\p_t \Om_0} + \nfy{\nabla_x \Om_0} + \eta_0 \nfy{D_0} ) \sum_{\substack{|\alpha|+|\beta|=N\\ 0\le |\beta|\le N}} \ep^{|\beta|+1} \nm{\nabla_\om^{\alpha} \nabla_x^\beta f}^2 \\\no
   & + \eta_0 \sum_{\substack{|\alpha|+|\beta|=N\\ 1\le |\beta|\le N}}
       \sum_{\substack{\beta_1+\beta_2=\beta\\ |\beta_1|\ge 1}} C_{\beta}^{\beta_1}
       \nfy{\nabla_x^{\beta_1} D_0} \nm{\nabla_\om^{\alpha+1} \nabla_x^{\beta_2} f} \nm{\nabla_\om^\alpha \nabla_x^\beta f} \ep^{|\beta|+1} \\\no
   & + \sum_{\substack{|\alpha|+|\beta|=N\\ |\alpha_1|\ge 1}}
       \sum_{\substack{\alpha_1+\alpha_2=\alpha\\\beta_1+\beta_2=\beta}} C_\alpha^{\alpha_1} C_{\beta}^{\beta_1}
       \nfyom{\nabla_\om^{\alpha_1} \nabla_x^{\beta_1} h_0}^2 \nm{\nabla_\om^{\alpha_2} \nabla_x^{\beta_2} f} \nm{\nabla_\om^\alpha \nabla_x^\beta f} \ep^{|\beta|+1} \\\no
   & + \sum_{\substack{|\beta|=N\\\beta_1+\beta_2=\beta}} C_{\beta}^{\beta_1}
        \nfyom{\nabla_x^{\beta_1} h_0}^2 \nm{\nabla_x^{\beta_2} f} \nm{\nabla_x^\beta f} \ep^{N+1}.
\end{align}
Denote
\begin{align*}
  &\mathcal{F}_N = \sum_{\substack{|\alpha|+|\beta|=N\\ 0\le |\beta|\le N}} \ep^{|\beta|+1} \nm{\nabla_\om^\alpha \nabla_x^\beta f}^2, \quad
   \mathcal{G}_N = \sum_{\substack{|\alpha|+|\beta|=N\\ 0\le |\beta|\le N}} \ep^{|\beta|} \nm{\nabla_\om^{\alpha+1} \nabla_x^\beta f}^2,
\\
  &\mathcal{H}_N = \sum_{\substack{|\alpha|+|\beta|=N\\ 0\le |\beta|\le N}} \ep^{|\beta|} \nm{\nabla_\om^{\alpha} \nabla_x^\beta h_1}^2,
\\
  & C_N= \sup \Big\{ \sum_{\substack{\mathclap{\beta_1+\beta_2=\beta}\\|\beta_1|\ge 1}} C_{\beta}^{\beta_1}
        \nfy{\nabla_x^{\beta_1} \Om_0}^2,\ \eta_0 \sum_{\substack{\mathclap{\beta_1+\beta_2=\beta}\\\mathclap{|\beta_1|\ge 1}}} C_{\beta}^{\beta_1}
       \nfy{\nabla_x^{\beta_1} D_0},\ \sum_{\substack{\mathclap{\alpha_1+\alpha_2=\alpha}\\\mathclap{\beta_1+\beta_2=\beta}}} C_\alpha^{\alpha_1} C_{\beta}^{\beta_1} \nfyom{\nabla_\om^{\alpha_1} \nabla_x^{\beta_1} h_0}^2  \Big\},
\end{align*}
then it follows from the above inequality that
\begin{align}\label{Esm:a priori-high order}
  \frac{1}{2} \frac{\rm d}{{\rm d} t} \mathcal{F}_N + \frac{1}{2} \mathcal{G}_N
 \ls & C_N (N\mathcal{F}_N + (N-1)\ep \mathcal{F}_{N-1} + \cdots + 2 \ep^{N-2} \mathcal{F}_2 + \ep^{N-1} \mathcal{F}_1 ) \\\no
       & + \mathcal{H}_N + \mathcal{F}_N + C_0 \mathcal{F}_N \\\no
 & + \eta_0 C_N \mathcal{F}_N^\frac{1}{2} (N \ep^\frac{1}{2} \mathcal{F}_N + (N-1)\ep \mathcal{F}_{N-1}^\frac{1}{2} + \cdots + \ep^\frac{N-1}{2} \mathcal{F}_1 + \ep^\frac{N}{2} \mathcal{F}_0 ) \\\no
 \le& C C_N (\mathcal{F}_N + \mathcal{F}_{N-1} + \cdots + \mathcal{F}_1 + \mathcal{F}_0) + C \mathcal{H}_N.
\end{align}


\end{remark}


\begin{remark}
  In fact, recall the definition
\begin{align*}
h_1=&-\frac{1}{M_0} \left[ \p_t f_1 + \om \cdot \nabla_x f_1 +\eta_0 \nabla_\om \cdot \left( P_{\om^\bot} D_{0} f_1 \right)  \right] \\
=&-\frac{1}{M_0} \left[ \p_t f_1 + \om \cdot \nabla_x f_1 -\eta_0 \left( (n-1) D_{0} \om f_1 - D_0 \nabla_\om f_1 \right)  \right],
\end{align*}
then by simple calculations, we can infer from Lemma \ref{lemm: Estimate for auxiliary} that
\begin{align*}
    \mathcal{H}_0(t) = & \nm{h_1}^2 \le C(\nm{\p_t f_1}, \nm{\nabla_x f_1}, \nm{\nabla_\om f_1}) \\
  \le & C (\|\rho_0\|_{L^\infty(0,T; H^m(\mathbb{T}))}, \|\Om_0\|_{L^\infty(0,T; H^m(\mathbb{T}))}),
\\
    \mathcal{H}_1(t) = & \ep \nm{\nabla_x h_1}^2 +\nm{\nabla_\om h_1}^2 \le C (\|\rho_0\|_{L^\infty(0,T; H^m(\mathbb{T}))}, \|\Om_0\|_{L^\infty(0,T; H^m(\mathbb{T}))}),
\\
    \mathcal{H}_2(t) = & \ep^2 \nm{\nabla_x^2 h_1}^2 + \ep \nm{\nabla_\om \nabla_x h_1}^2 + \nm{\nabla_\om^2 h_1}^2\\
    \le & C (\|\rho_0\|_{L^\infty(0,T; H^m(\mathbb{T}))}, \|\Om_0\|_{L^\infty(0,T; H^m(\mathbb{T}))}),
\end{align*}
for any $t\in[0,T]$ and $m>13$.

As a consequence, by choosing a new constant $\hat{C} \deq \sup_{t\in[0,T]}\{ C_0, C_1, C_2, \mathcal{H} \}$, we can fet a refined formulation of \eqref{eq: a priori},
\begin{align}
  \frac{\d}{\d t} \mathcal{E} + \mathcal{G} \le C \hat{C} (\mathcal{E} + 1),
\end{align}
which can lead to the same result. Here we keep the original formulation with the term $\mathcal{H}$ in order to make the contributions of $f_1$ more clear.
\end{remark}

\section{Completion of the Proof of the Main Theorem} 
\label{sec:completion_of_main_thm}

\subsection{Local Existence of the Remainder Equation} 
\label{sub:local_exist_of_remainder_eq}


Based on the a priori estimates (Lemma \ref{lemm: a priori}), we can get the local existence for the remainder equation \eqref{Eq: Target} by a standard iteration scheme.
\begin{lemma}\label{lemm: local exist}
  Given the initial datum $f(0,x,\om) \in H^2(\mathbb{T} \times \mathbb{S}^{n-1})$, then there exist $T_*>0$ such that the remainder equation \eqref{Eq: Target} admits a unique solution $f(t,x,\om)$ on $C([0,T_*); H^2(\mathbb{T} \times \mathbb{S}^{n-1}))$, and moreover, there exists a constant $E>0$,
  \begin{align*}
    \sup_{t\in[0,T_*]} \mathcal{E}(f(t)) \le E,
  \end{align*}
  provided that the initial datum satisfies $\mathcal{E}(f(0)) \le \frac{E}{2}.$
\end{lemma}

Observing that the equation \eqref{Eq: Target} is linear with respect to $f$, proving the lemma can be established by a standard iteration scheme and a straightforward compactness justification, hence we will only sketch the proof. Here we consider the following iteration scheme:
\begin{align}\label{Eq:iteration} \tag{$P_n$}
  \left\{
    \begin{array}{l}
      \p_t (\ep f^{n+1}) + \om \cdot \nabla_x (\ep f^{n+1}) +\eta_0\, \nabla_\om \cdot \left( P_{\om^\bot} D_{0} (\ep f^{n+1}) \right) + \mathcal{L}_0 f^{n+1}
      = h_1 + \ep f^n h_0, \\[0.5em]
      f^{n+1}(0,x,\om) =f(0,x,\om),
    \end{array}
  \right.
\end{align}
starting with $f^0(0,x,\om) =f(0,x,\om)$.

We remark that under the hypothesis of Theorem \ref{Thm: main}, there exists $\ep_0 >0$ such that,
\begin{align*}
  \mathcal{E}(f^{n+1}(0)) = \mathcal{F}^{n+1}_0(0)+\mathcal{F}^{n+1}_1(0) +\mathcal{F}^{n+1}_2(0) \le \frac{E}{2}
\end{align*}
holds for the above $E>0$ and all $\ep \in (0, \ep_0) $.

To complete the proof of Lemma \ref{lemm: local exist}, it suffices to get a uniform (in $n$) estimate for $\mathcal{E}^{n+1}(t) $.

\begin{lemma}\label{lemm: uni-bound-4-iter}
  There exists $T_*>0$, such that if $\sup_{t\in[0,T_*]}\mathcal{E}(f^n(t)) \le E$, then
  \begin{align*}
    \sup_{t\in[0,T_*]}\mathcal{E}(f^{n+1}(t)) \le E.
  \end{align*}
\end{lemma}
\proof Following exactly the same lines as the proof of the a priori estimates, we can get finally
\begin{align*}
  &\frac{\d}{\d t} \mathcal{F}^{n+1}_0 + \mathcal{G}^{n+1}_0 \ls C_0 (\mathcal{F}^{n+1}_0 + \mathcal{F}^n_0) + \mathcal{H}_0, \\
  &\frac{\d}{\d t} \mathcal{F}^{n+1}_1 + \mathcal{G}^{n+1}_1 \ls C_1 (\mathcal{F}^{n+1}_1 + \mathcal{F}^n_1 + \mathcal{F}^n_0) + \mathcal{H}_1, \\
  &\frac{\d}{\d t} \mathcal{F}^{n+1}_2 + \mathcal{G}^{n+1}_2 \ls C_2 (\mathcal{F}^{n+1}_2 + \mathcal{F}^n_2 + \mathcal{F}^n_1 + \mathcal{F}^n_0) + \mathcal{H}_2,
\end{align*}
hence a constant $\widetilde{C}>0$ depending on $\sup_{t\in[0,T_\ep]} \{C_0, C_1, C_2\}$ exists such that,
\begin{align}
  \frac{\d}{\d t} \mathcal{E}(f^{n+1}) + \mathcal{G}(f^{n+1})
  \le \widetilde{C} (\mathcal{E}(f^{n+1}) + \mathcal{E}(f^n) + \mathcal{H}).
\end{align}

Integrating from 0 to $t$ with respect to the time variable, we have
\begin{equation}\nonumber
\begin{aligned}
    &\mathcal{E}(f^{n+1}(t)) + \int_0^t \mathcal{G}(f^{n+1}(s)) \d s\\
  \le & \mathcal{E}(f^{n+1}(0)) + \widetilde{C} t (\sup_{s\in [0,t]} \mathcal{E}(f^{n+1}(s)) + \sup_{s\in [0,t]} \mathcal{E}(f^n(s)) + \sup_{s\in [0,t]} \mathcal{H}(s)).
\end{aligned}
\end{equation}

Take sufficiently small $T_*$ such that,
\begin{align*}
  \widetilde{C} T_* \le \frac{1}{6}, \quad
  \widetilde{C} T_* \sup_{s\in [0,t]} \mathcal{H}(s) \le \frac{1}{6} E,
\end{align*}
which together with the assumptions $\sup_{t\in[0,T_*]}\mathcal{E}(f^n(t)) \le E$ and $\mathcal{E}(f^{n+1}(0)) \le \frac{1}{2} E$ enables us to get the desired result
\begin{align}
  \sup_{t\in[0,T_*]}\mathcal{E}(f^{n+1}(t)) \le E.
\end{align}
\qed

Then we are left to prove the convergence. Set $w^n= f^{n+1} - f^n$, for which the iteration scheme \eqref{Eq:iteration} leads to
\begin{align}
  \left\{
    \begin{array}{l}
      \p_t (\ep w^n) + \om \cdot \nabla_x (\ep w^n) +\eta_0\, \nabla_\om \cdot \left( P_{\om^\bot} D_{0} (\ep w^n) \right) + \mathcal{L}_0 w^n
      = \ep w^{n-1} h_0, \\[0.5em]
      w^n(0,x,\om) =0.
    \end{array}
  \right.
\end{align}
The same computations as used for \eqref{Eq:iteration} give that
\begin{align*}
    \mathcal{E}(w^n(t)) + \int_0^t \mathcal{G}(w^n(s)) \d s
  \le \widetilde{C}T_* (\sup_{s\in [0,T_*]} \mathcal{E}(w^n(s)) + \sup_{s\in [0,T_*]} \mathcal{E}(w^{n-1}(s))).
\end{align*}
Take $T_*$ small enough so that $\widetilde{C} T_* \le \frac{1}{3} $, then we obtain
\begin{align*}
    \sup_{t\in [0,T_*]} \mathcal{E}(f^{n+1}(t) - f^n(t)) \le \frac{1}{2} \sup_{t\in [0,T_*]} \mathcal{E}(f^n(t) - f^{n-1}(t)).
\end{align*}
Consequently, it can be proved easily the sequence $\{f^n\}$ is convergent in the space $C(0,T_*; $ $H^2(\mathbb{T} \times \mathbb{S}^{n-1}))$ and so is its limit $f$. This completes the proof of Lemma \ref{lemm: local exist}.



\subsection{Completion for the Main Theorem} 
\label{sub:justify_main_thm}


We are now in a position to complete the proof of our main Theorem \ref{Thm: main}. Lemma \ref{lemm: local exist} shows that there exists a maximal time $T_\ep>0$ such that the solution $f_2^\ep$ to the remainder equation \eqref{Eq: Target} satisfies the a priori estimate stated in lemma \ref{lemm: a priori},
\begin{align}
\frac{\d}{\d t} \mathcal{E} + \mathcal{G}
\le \widetilde{C} (\mathcal{E} + \mathcal{H}),
\end{align}
for any $\ep \in (0,\ep_0)$ and $t\in[0,T_\ep] $.

Observing the hypothesis of Theorem \ref{Thm: main} $\|f_2^\ep\|_{H^2_{x,\om}} (0) \le C$ implies that $\mathcal{E}(0) \le E$. Define
\begin{align*}
  E_2 = e^{\widetilde{C} T}(E(0) + \widetilde{C} T \sup_{t \in [0,T]} \mathcal{H}(t)), \\
  T_0= \sup \{ t \in [0,T_\ep],\ \mathcal{E}(t) \le E_2 \}.
\end{align*}

We claim that $T_\ep \ge T$. Indeed, applying the Gr\"onwall inequality, we infer from the above a priori estimate, for $t\in[0,T_\ep] $,
\begin{align}
  \mathcal{E}(t) \le e^{\widetilde{C} t}(E(0) + \widetilde{C} t \sup_{s \in [0,t]} \mathcal{H}(s)).
\end{align}
If $T_\ep < T$, then it follows that
\begin{align*}
  \mathcal{E}(t) \le E_2,
\end{align*}
which in turn yields that $T_\ep = T_0$, and furthermore, the a priori estimate are adapted to the time $t= T_\ep$ so that the solution may be continued beyond $T_\ep$. This stands in contradiction to the maximal property of $T_\ep$. Thus $T_\ep \ge T$, and $\mathcal{E}(t) \le E_2$ holds for $t\in [0,T]$.

Recalling the definition of $\mathcal{E}$
\begin{align*}
   & \mathcal{E} = \mathcal{F}_0 + \mathcal{F}_1 + \mathcal{F}_2 \\\no
  =& \ep \nm{f}^2 +(\ep^2 \nm{\nabla_x f}^2 + \ep \nm{\nabla_\om f}^2) + (\ep^3 \nm{\nabla_x^2 f}^2 + \ep^2 \nm{\nabla_\om \nabla_x f}^2 + \ep \nm{\nabla_\om^2 f}^2)
\end{align*}
and the equivalence between the two norms $\nm{\cdot}$ and $\|\cdot\|$, we actually have proved the result \eqref{Eq: conclusion},
\begin{align}
  \ep^\frac{1}{2} \|f_2^\ep\|_{L^2_{x,\om}} + \ep \|f_2^\ep\|_{H^1_{x,\om}}+ \ep^\frac{3}{2} \|f_2^\ep\|_{H^2_{x,\om}} \le C,
\end{align}
with $C$ is independent of $\ep \in (0,\ep_0) $ and $t\in[0,T] $.





\end{document}